\documentclass[11pt]{article}              % for a regular run

\usepackage{amsmath,amstext,amsfonts,amssymb}%amsthm
\usepackage{graphicx,color}
\usepackage{personal2}

\newcommand{\la}{\lambda}
\newcommand{\eps}{\varepsilon}

\newcommand{\al}{\alpha}

\newcommand{\sig}{\sigma}
\newcommand{\del}{\delta}

\newcommand{\La}{\mathnormal{\Lambda}}

\newcommand{\Ph}{\mathnormal{\Phi}}

\newcommand{\N}{{\mathbb N}}

\newcommand{\R}{{\mathbb R}}

\newcommand{\Z}{{\mathbb Z}}

\newcommand{\E}{{\mathbb E}}
\newcommand{\PP}{{\mathbb P}}

\newcommand{\calF}{{\cal F}}

\newcommand{\calL}{{\cal L}}

\newcommand{\calT}{{\cal T}}
\newcommand{\calU}{{\cal U}}

\newcommand{\skp}{\vspace{\baselineskip}}

\newcommand{\To}{\Rightarrow}

\newcommand{\iy}{\infty}
\newcommand{\noi}{\noindent}
\newcommand{\qed}{\hfill$\Box$}
\newcommand{\sDplus}{\sD^+}

\providecommand{\KEYWORDS}[1]{\textbf{\textit{Index terms---}} #1}
\providecommand{\MSCCLASS}[1]{\textbf{\textit{MSC2010 subject
      classification---}} #1}
\providecommand{\ORMSCLASS}[1]{\textbf{\textit{OR/MS subject classification---}} #1}

\title{Asymptotically optimal appointment schedules with customer no-shows}
\author{Mor Armony\thanks{Stern School of Business,
New York University {\tt marmony@stern.nyu.edu}}
\and
Rami Atar\thanks{Viterbi Faculty of Electrical Engineering,
Technion--Israel Institute of Technology {\tt atar@ee.technion.ac.il}}
\and
Harsha Honnappa\thanks{School of Industrial Engineering,
Purdue University {\tt honnappa@purdue.edu}
\newline
\newline
 \MSCCLASS{Primary: 68M20, 90B36; Secondary: 60K25, 90C15}
 \newline
 \ORMSCLASS{Primary: queues: optimization; Secondary:
   production/scheduling : sequencing stochastic; programming:
   stochastic}
}
}

\date{\today}

\begin{document}

\maketitle

\begin{abstract}
We consider the problem of scheduling appointments for a
finite customer population to a service facility with customer
no-shows, to minimize the sum of customer waiting time and server overtime
costs. Since appointments need to be scheduled ahead of time we refer to this problem as an optimization problem rather than a dynamic control one. We study this optimization problem in fluid and diffusion scales and identify asymptotically optimal schedules in both scales. In fluid scale, we show that
it is optimal to schedule appointments so that the
system is in critical load; thus
heavy-traffic conditions are obtained as a result of
{\it optimization} rather than as an {\it assumption}. In diffusion scale, we solve this optimization problem in the large horizon limit.
Our explicit stationary solution of the corresponding Brownian Optimization Problem translates the customer-delay versus server-overtime tradeoff to a tradeoff between the state of a
reflected Brownian motion in the half-line and its local time at zero.
Motivated by work on competitive ratios,
we also consider a reference model in which an oracle provides
the decision maker with the complete randomness information.
The difference between the values of the scheduling problem for the two models,
to which we refer as the {\it stochasticity gap} (SG),
quantifies the degree to which it is harder to design a schedule under uncertainty
than when the stochastic primitives
(i.e., the no-shows and service times) are known in advance.
In the fluid scale, the SG converges to zero,
but in the diffusion scale it converges to a positive constant that we compute.
\end{abstract}

\section{Introduction}
We study the problem of determining an optimal appointment schedule for a
finite number of customers at a service system that only accepts arrivals in a finite time horizon, but renders
service to all arriving customers. Broadly, the objective here is to assign deterministic arrival epochs to a finite
population such that the server is optimally utilized while the
cumulative delay experienced by the customers is minimized. The optimization problem is stochastic in nature due not only to the
randomness in service times, but also to the fact that some arrivals do not
show up. As opposed to typical stochastic control problems that appear in the literature, the schedule needs to be determined {\em offline}, ahead of time, with no access to the realization of the stochastic primitives over time.

Our research is motivated by systems such as outpatient clinics
that render service to a finite number of patients during a working
day (7AM to 4PM, for instance). No patients are accepted for service if they arrive after
the end of the horizon, but all patients that do arrive are rendered
service. No-shows in outpatient care is a problem most clinics
struggle with regularly. According to \cite{CaVeRo2006} no-show rates may be up to 60\%, depending on the clinic specific characteristics. Patient overbooking has been proposed as an effective strategy to handle clinic under-utilization resulting from patient no-show (\cite{LaLa2012}). At the same time overbooking may lead to clinic overcrowding, that will intensify patient waiting time and doctor overtime. Our paper aims at determining effective appointment schedules that will minimize these wait- and over- times.

We model the service system as a single server queue with an infinite
buffer. We assume that the service times are
generally distributed IID non-negative random variables (RVs) with a finite second moment, and that the server
is non-idling and operates according to a first-come-first-served discipline.
The entire finite population needs to be provided with appointments, and these appointments are allowed to span the entire time horizon, but not beyond.
The actual arrival process is thinned with probability $p$ due to no-shows from this deterministic
appointment schedule.

Our optimization problem's goal is to determine an appointment schedule to minimize the objective of a weighted sum of the
expected cumulative wait time of all the customers that arrive in the
(finite) arrival horizon and the
expected overage time, defined as the amount of time it takes for
the server to clear out the backlog after the end of the horizon.
The appointment scheduling problem is reducible
to a stochastic bin-packing problem making it
NP-hard, and is generally solved using various heuristics
(see \cite{GuDe2008}).

Here, we introduce two large population limit regimes
that lead to
simpler optimization problems and yield exact solutions.
We operate in a large population limit framework that reveals the fundamental
  complicating factors in the optimization problem. Our scaling
  regimes let the population size tend to infinity while
  simultaneously accelerating the service rate in proportion to the
  population size. In the fluid regime, the service time cost is scaled by the inverse of the population size.  Customers are assigned arrival epochs according to a sequence of
arbitrary schedules. In the limit, the fluid regime washes out
the stochastic variation and captures the `mean' or first order
effects in the queue performance.  We posit a
variational fluid
optimization problem (FOP) and solve it
in Proposition~\ref{thm:fop}, under an `overload' condition that the aggregate available fluid
service is less than the expected aggregate fluid arrivals, and thus the overtime cost is non-zero. The optimal cumulative fluid schedule function matches the cumulative service
completions in the arrival horizon and schedules the remaining fluid at
the end of the horizon. This result shows that the
heavy-traffic critical-load condition emerges as a consequence of optimization as opposed to a postulated assumption. In
Theorems~\ref{prop:fluid-cost-lower-bound} and~\ref{thm:asymptotic-achievability} we prove that the value of
the fluid optimization problem is asymptotically achievable by a
carefully constructed sequence of simple finite population schedules.

While the FOP results in a simple and intuitive schedule, it only considers first-order effects and washes away stochastic fluctuations. In actuality, we expect the inherent stochasticity of the system to
have a significant impact on system performance and thus on the design of the schedule. Indeed,
one might ask, whether considering second-order terms will shed more light on the optimal appointment schedule. Specifically, it is of interest to see how fluctuations of order of the square root of the population size about the fluid solution impact the schedule. To formalize this question we posit
a Brownian optimization problem (BOP), assuming
the same overload condition as in the fluid scale.
The BOP is stated in terms of an equation driven by a Brownian motion and a
control, obtained as a formal diffusion limit of the queue length in heavy traffic.
This is in line with the fact that the fluid optimal schedule
enforces criticality within the appointment horizon. The BOP is not a dynamic control problem, but one where the control trajectory has to be planned ahead of time zero. This makes insights, structure, and tools from dynamic programming irrelevant.
While we are unable to obtain an explicit solution of the BOP when set on a fixed time horizon,
we do derive an explicit solution to it in the large time horizon limit.
In Proposition~\ref{prop4} we prove that the value of the BOP in the
large horizon limit can be achieved by an RBM that has a constant negative
drift. Equivalently, in this limit, the BOP is
solved by a stationary RBM. We identify the optimal drift coefficient in Lemma~\ref{lem4--}.
In Theorems~\ref{prop5} and~\ref{prop:diff-opt-sched} we prove that the value of
the BOP is also asymptotically achievable by a
carefully constructed sequence of finite population schedules.

There is, of course, a `price' to be paid for having to schedule traffic at time
zero without any stochastic information revealed ahead of time. We quantify this by introducing the notion of a {\it stochasticity gap} (SG), defined as the
difference between the appropriately scaled finite population value
and the value of the `complete information' (CI) problem. In the CI
problem an oracle reveals all stochastic primitives (or future events)
  to the optimizer at time zero. The CI problem is not completely trivial, but much easier than
  the original one. The CI and the fluid optimal schedules are similar in that they both schedule appointments such that customers arrive at (precisely or approximately, respectively) the time when they are ready to be served, and
any excess jobs are scheduled at the end of the horizon.
Indeed, Proposition~\ref{thm:sg-fluid} shows that in the fluid scale the
asymptotic SG is zero. On the other hand, in diffusion scaling we calculate the SG and show
in Proposition~\ref{thm:sg-diffusion} that it
is strictly positive.
(The term SG has been used before in the robust optimization literature
in relation to performance with and without information revealed by an oracle,
but in a way that differs from our use of it, see eg.\ \cite{BG}.)

To summarize, our main contributions are as follows:
\begin{enumerate}
\item[(i)] Under an overload condition, we identify explicitly computable
  asymptotically optimal
  schedules in the fluid and diffusion regimes. These constitute the
  first analytical results for the appointment scheduling problem with
  no-shows in the large population limit.

\item[(ii)] Most literature on queueing control problems in heavy traffic
  focuses on dynamic control. Our
  results differ from this literature in two fundamental
  ways. First, the queueing model problem is formulated as an optimization
  rather than a dynamic control problem.
  Second, a critical load condition, which sets the ground for a heavy traffic
  analysis at the diffusion scale,
  emerges because it is optimal to operate at criticality, rather than as an assumption.

\item[(iii)] The essence of the optimization problem
  is that it must be carried out without the randomness being revealed.
  We analyze the SG as a means of quantifying the cost
  associated with uncertainty. To the best of our knowledge, this is the first
  work to study the SG at any scaling limit.
\end{enumerate}

The rest of the paper is organized as follows. We conclude this section with a brief overview of the relevant literature and notational convention. \S \ref{sec:model} provides the problem formulation and a summary of the main results. In \S\S \ref{sec:problem} and~\ref{sec:diff} we solve the problem under fluid and diffusion scaling, respectively. We compute the corresponding SG in these sections as well. We conclude with final remarks and a discussion of future research directions in \S \ref{sec5}.

\skp

We now review some of the relevant results in the field.
There is a vast literature on
appointment scheduling in healthcare which we will not attempt to
summarize here; we direct the reader to the comprehensive reviews in
\cite{CaVe2003,GuDe2008,Ha2012}. We note two that are particularly
relevant to our study. A scheduling problem close to ours has been
studied in \cite{ZaPi2014}. The problem they consider is that of
determining an optimal schedule for heterogeneous patients in the
presence of no shows with the objective of minimizing waiting cost
plus idling and overtime cost. Our model may be considered a special
case of their model in that our patients are homogeneous, and our
idling cost is zero. However, their model considers fixed appointment
slots, while ours allows for appointment times to be a result of the
optimization problem. For a finite and fixed patient population, the
authors characterize structural properties of an optimal
schedule. Explicit solutions are given for some special cases and are
studied numerically for the more general problem. The authors observe
that optimal solutions tend to front load (more overbooking towards
the beginning of the day). The numerical solutions of \cite{ZaPi2014} also show that some appointment slots should be overbooked but not all, with only up to a couple of patients scheduled per slot.

A second paper that is relevant is \cite{HaMe2008}, which considers a similar problem definition. However, the model there assumes that service times are exponentially distributed and there is no fixed horizon in which the finite number of arrivals must be scheduled. The cost function again trades off the expected cumulative waiting time of the customers that show-up against the expected ‘server time’ beyond the last scheduled arrival epoch. The authors provide extensive numerical analysis of the finite-population scheduling problem. In particular, they numerically compute the optimal schedule which shows that overbooking is possible for the first few and last few arrival epochs, and arrivals ‘in the middle’ are almost uniformly spread out. Furthermore, they also contrast the value of their problem against that of an oracle problem akin to our CI problem, and note the fact that value of the latter is significantly lower.

 Also relevant is \cite{BeJo2010} who consider finite population scheduled arrivals models with both no-shows and tardiness. Under the assumption of exponentially distributed service times, the authors derive exact expressions for various performance metrics.

The second relevant stream of literature is on the asymptotics of
scheduled arrivals (\cite{ArGl2012}) and transitory queueing
(\cite{HoJaWa2015,HoJaWa2016}). This stream considers limits of
processes that are generated from a queue with an arrival process that
is originated from a finite population. Closest to our model is \S 4.3
of \cite{HoJaWa2016}, in which the authors consider scheduled arrivals with epoch uncertainty. Indeed, the asymptotic scaling and limiting regimes are similar to ours. The two main differences are that (i) the appointment times are given and are assumed a-priori to be equally spaced, and (ii) all customers are assumed to show up but they may be non-punctual.

\paragraph{Notation}
Let $\sD[0,\iy)$ be the space of functions $f: [0,\iy) \to \bbR$ that are right continuous
with left limits, equipped with the Skorohod $J_1$ topology.
Let $\sDplus[0,\iy) \subset \sD[0,\iy)$ be the subset of
non-negative non-decreasing functions.
For a sequence $\{X_n\}$, $X$,
of RVs, $X_n \Rightarrow X$ as $n \to \infty$ denotes convergence in law.
For a sequence $\{X_n\}$, $X$, of stochastic processes with sample paths
in $\sD[0,\iy)$, $X_n\To X$ denotes convergence in law in the $J_1$
topology. In this paper, all statements involving convergence of processes $X_n\To X$
will be to processes $X$ that have a.s.\ continuous sample paths,
thus these convergences can equivalently be understood as convergence in law
in the {\it uniformly on compacts} (u.o.c.) topology.
Let $(\cdot)^+ := \max\{\cdot,0\}$.
For an event $A$, $\mathbf 1_{A}$ is the corresponding indicator function.
For $f:\R_+\to\R$, $T>0$, $\del>0$ we denote $\|f\|_T=\sup_{t\in[0,T]}|f(t)|$. A one-dimensional Brownian motion (BM) with drift $m$ and diffusion coefficient $\sigma$,
starting from zero, is referred to as an $(m,\sigma)$-BM.
The letter $c$ denotes a positive constant whose value is immaterial,
and may change from line to line.

\section{Problem setting} \label{sec:model}
\subsection{Model}
Consider a single server queue with an
infinite waiting room. A finite number of jobs arrive at
the queue over a finite time horizon, and are served
on a first-come-first-served (FCFS) schedule. Jobs are given appointments at fixed times
during the day (not necessarily uniformly spaced), and we assume that
jobs that do turn up do so precisely at the appointment time; that is, we assume punctual arrivals but allow no-shows. We also assume that the service times are independent and identically  distributed (IID) with finite second moments.

The RVs and stochastic processes
are defined on a probability space $(\Omega,\calF,\PP)$, and the symbol
$\E$ denotes expectation with respect to $\PP$. The number of requested appointments (or population size) is denoted as $N$.
Let $H > 0$ denote the operating time horizon.
A {\it schedule} is any deterministic, non-decreasing sequence
$\{T_i, ~i = 1,\ldots, N\}$ taking values in $[0,H]$.
It represents scheduled arrival epochs. The collection of all schedules is denoted by $\calT$.
We denote by $E(t)$ the cumulative number of scheduled arrivals by
time $t$, that is, $E(t) = \sum_{i=1}^{N} \mathbf 1_{\{T_i \leq t\}}$. The function
$E$ is referred to as the {\it scheduling} function.
Let $\{\xi_i,~i\in\N\}$ be IID Bernoulli RVs with mean
$p\in(0,1]$. They are used as model for actual arrivals, namely
$\xi_i=1$ if and only if the $i$th scheduled job shows up.
With these elements, the {\it cumulative arrival}
process, $A \in \sDplus[0,H]$, is given by
\(
A(t)=\sum_{i=1}^{N} \xi_i\mathbf 1_{\{T_i \leq t\}}.
\)
Note that, with $\Xi(k) = \sum_{i=1}^k \xi_i$, $k\in\Z_+$,
one can express the above relation as
\begin{equation}
  \label{eq:arrival}
  A(t) = \sum_{i=1}^{E(t)}\xi_i= \Xi \circ E(t).
\end{equation}

Let $\{\nu_i\}$ be an IID sequence of non-negative RVs with mean $\mu^{-1}$
and squared coefficient of variation $C_S^2\in(0,\iy)$.
Assume that this sequence and the sequence $\{\xi_i\}$ are mutually independent.
These RVs are used to model service times in the following way.
Let $\{\nu_i, ~i\in\N\}$ be the service time of the $i$th served job.
Let
$S(t) = \max \{m | \sum_{i=1}^m\nu_i \leq t\}$.
Then $S(t)$ is the cumulative number of service completions
by the time the server is busy for $t$ units of time.

Let $Q$ denote the {\it number-in-system} process. Then the {\it cumulative busyness}
process is given by $B(\cdot) = \int_0^\cdot \mathbf 1_{\{Q(s) > 0\}} ds$.
A simple balance equation for $Q$ is
\begin{equation}
  \label{eq:queue-length}
  Q = A - S \circ B.
\end{equation}

\subsection{Cost and optimization problem}

Two primary performance measures of interest to us are the {\it overall
waiting time}, or {\it makespan}, and the {\it overage time}.
The former is defined as
\begin{equation}
  \label{eq:waiting-time}
  W = \int_0^\infty (Q(s) - 1)^+ ds,
\end{equation}
and represents the sum, over all arriving jobs, of the job's waiting time in the queue
(not counting the time of service).
The overage time is the amount of time after the end of the horizon $[0,H]$ it takes
to complete the last arrival in $[0,H]$. If we denote by $\tau$
the time when the last arrival in $[0,H]$ departs from the server,
then the overage time is given by
\begin{equation}\label{eq:overage}
O=(\tau-H)^+.
\end{equation}
Note that $\tau$ can be expressed in terms of the processes introduced earlier as
\begin{equation}
  \label{50}
\tau=\inf\{t:S\circ B(t)\ge \Xi(N)\},
\end{equation}
for $\Xi(N)$ is the total number of arrivals, and $S\circ B(t)$ is the number
of departures by time $t$.

The goal of the system operator is to schedule the $N$ jobs so as to
minimize a combination of the expected overall waiting time experienced by the jobs (makespan)
and the expected overage time. Thus we consider
the {\it finite population optimal-schedule problem (FPOP)} defined by considering the cost
\begin{equation} \label{eq:FPOP}
J(\{T_i\}) = c_w \bbE[W] + c_o \bbE [O],\qquad \{T_i\}\in\calT,
\end{equation}
and the value
\begin{equation}\label{51}
V=\inf_{\{T_i\}\in\calT} J(\{T_i\}).
\end{equation}
Here, $c_w$ and $c_o$
are non-negative constants.
In general, solving the FPOP is formidable. Rather than
solve it directly, we solve  fluid and diffusion scale problems and
prove the existence of asymptotically optimal finite population
schedules that approach the fluid- and diffusion-optimal solutions in
the large population limit.

The main objective of this work is to study the asymptotics
of problem \eqref{eq:FPOP}--\eqref{51}. Before we turn to it
we comment on another, related problem setting.
Suppose that all the information on the stochastic data is known
to the decision maker when he selects the schedule at time zero.
In this version of the problem, that we refer to as the {\it complete information} (CI)
problem, the selection of schedule $\{T_i\}$ may depend on
the stochastic data $(\{\xi_i\},\{\nu_i\})$.
We do not consider it as a practically motivated setting by itself, because in applications
we have in mind these stochastic ingredients are not known in advance.
However, it is useful to regard it as a reference model, and to
relate it to our main problem \eqref{eq:FPOP}--\eqref{51}.

For a precise
formulation, let $\Sigma$ denote the sigma-field generated by the collection
$(\{\xi_i\},\{\nu_i\})$. An $\Sigma$-measurable RV $\{\Theta_i\}$
taking values in $\calT$ is called a {\it CI-schedule}.
Let $\calT^{\rm CI}$ denote the collection of all CI-schedules
(note that schedules in $\calT^{\rm CI}$ are not allowed to change the {\em order} of the scheduled jobs).
Then, analogously to \eqref{eq:FPOP}, we let
\begin{equation}
  \label{52}
  J^{\rm CI}(\{\Theta_i\})=c_w\E[W]+c_o\E[O],\qquad \{\Theta_i\}\in\calT^{\rm CI},
\end{equation}
\begin{equation}
  \label{53}
  V^{\rm CI}=\inf_{\{\Theta_i\}\in\calT^{\rm CI}}J^{\rm CI}(\{\Theta_i\}),
\end{equation}
where $W$ and $O$ correspond to the selection $\{\Th_i\}$.
Clearly, one always has $V\ge V^{\rm CI}$. We refer to the difference
\[
\gamma=V-V^{\rm CI}
\]
as the {\it stochasticity gap}, as it quantifies the gap between performance
with and without knowing the stochastic ingredients.

Problem \eqref{52}--\eqref{53} is much easier than our main problem,
and is in fact fully solvable.
We devote \S \ref{sec24} to present its solution.

\subsection{Large-population asymptotic framework}\label{subsec_asym}
Because the problem \eqref{50}--\eqref{51} is prohibitively difficult to solve exactly, we instead take a large-population asymptotic approach where we consider a sequence of systems where the number of scheduled jobs grows large and the cost is scaled to make the problem tractable. Specifically, we consider a sequence of systems
indexed by $n\in\N$, where the population size $N_n$ satisfies $N_n=\left\lceil\al n\right\rceil$, where $\al>0$ is a fixed parameter.
A schedule in the $n$th system is a non-decreasing sequence $\{T_{i,n},
~i = 1,\ldots, N_n\}$ taking values in $[0,H]$, with the corresponding
scheduling {function} defined as  $E_n(t) = \sum_{i=1}^{N_n} \mathbf 1_{\{T_{i,n} \leq t\}}$, and $A_n(t) = \sum_{i=1}^{E_n(t)}\xi_i= \Xi \circ E_n(t)$.
The collection of all schedules for the $n$th
system is denoted by $\calT_n$.
Let $S_n(t):=S(nt)=\max \{m | \sum_{i=1}^m\nu_i \leq nt\}= \max \{m | \sum_{i=1}^m\nu_{i,n} \leq t\}$, where $\nu_{i,n}=n^{-1}\nu_i$.
In parallel to (\ref{eq:queue-length}), the resulting number-in-system process for the $n$th system may be expressed as
\begin{equation}
  \label{eq:queue-length-n}
  Q_n = A_n - S_n \circ B_n,
\end{equation}
with $B_n(\cdot) = \int_0^\cdot \mathbf 1_{\{Q_n(s) > 0\}} ds$.

Under the large-population scaling, population is assumed to grow linearly with $n$. One may interpret the scaling of the various processes a couple of different ways. Under one interpretation time is scaled by $n$, the scheduling time horizon becomes $[0,nH]$, and service times are of order ${\mathcal O}(1)$. A second interpretation is that the scheduling time horizon remains as $[0,H]$, but the service times are scaled by $1/n$. That is, when the population grows, the server speeds up at a rate that is proportional to the population size. While these two interpretations are mathematically equivalent, we will provide intuition throughput the paper that is consistent with the second interpretation.

Paralleling (\ref{eq:waiting-time}), (\ref{eq:overage}), and (\ref{50}), we have that the makespan, the overage time, and the departure time of the last arrival, are respectively defined as:
\begin{equation}
  \label{eq:waiting-time-n}
  W_n = \int_0^\infty (Q_n(s) - 1)^+ ds,
\end{equation}
\begin{equation}
\label{eq:overage-n}
O_n=(\tau_n-H)^+,
\end{equation}
and
\begin{equation}
  \label{50-n}
\tau_n=\inf\{t:S_n\circ B_n(t)\ge \Xi(N_n)\}.
\end{equation}

Paralleling (\ref{eq:FPOP}) and (\ref{51}),
the {\it large population optimal-schedule problem (LPOP)} defined by considering the cost
\begin{equation} \label{eq:FPOP-n}
J_n(\{T_{i,n}\}) = c_{w,n} \bbE[W_n] + c_{o,n} \bbE [O_n],\qquad \{T_{i,n}\}\in\calT_n,
\end{equation}
where $c_{w,n}$ and $c_{o,n}$ are appropriately scaled constants,
and the value
\begin{equation}\label{51-n}
V_n=\inf_{\{T_{i,n}\}\in\calT_n} J_n(\{T_{i,n}\}).
\end{equation}

Similarly, for the complete information case we have
\begin{equation}
  \label{52n}
  J_n^{\rm CI}(\{\Theta_{i,n}\})=c_{w,n}\E[W_n]+c_{o,n}\E[O_n],\qquad \{\Theta_{i,n}\}\in\calT^{\rm CI},
\end{equation}
and
\begin{equation}
  \label{53n}
  V^{\rm CI}_n=\inf_{\{\Theta_{i,n}\}\in{\calT}_n^{\rm CI}}J_n^{\rm CI}(\{\Theta_{i,n}\}),
\end{equation}
Finally, let
\[
\gamma_n=V_n-V_n^{\rm CI}.
\]

To see what is the appropriate scaling for the cost coefficients of $J_n(\cdot)$, note that the leading (first-order) term in the expression for the number of jobs that the server can handle in the interval $[0,H]$ is $n\mu H$ (recall that we assume that in the $n$th system the server works at a rate $n\mu$). Similarly, the leading term in the total number of jobs that arrive into the system in $[0,H]$ is $p\alpha n$. To capture the case where the server incurs a non-negligible overage cost, we assume that the system is {\em overloaded}. That is, we assume that $p\alpha n > n\mu H$, or equivalently,
\begin{equation}\label{eq:overload}
p\alpha  > \mu H.
\end{equation}
Thus, regardless of the schedule, at time $H$ the number of jobs
present in the system is of order $n$. Since the service rate is also
of order $n$, it takes a constant (order $1$) time to handle these
jobs. That is, the overage time is ${\mathcal O}(1)$. Along the same
lines, notice that the number of arriving jobs is of order $n$ and
their individual waiting time is of order 1. Thus the total waiting time (makespan) is ${\mathcal O}(n)$. This suggests that in order to get a meaningful cost function $J_n(\cdot)$, the cost parameter $c_{w,n}$ should be scaled by $n^{-1}$ and the cost parameter $c_{o,n}$ should remain a constant. Thus, we assume for the rest of the paper that
\begin{equation}\label{eq:scaled_c}
c_{w,n}=n^{-1} c_w,~~c_{o,n}=c_o.
\end{equation}

We study this asymptotic problem under two scalings. The first is a fluid scaling (see \S \ref{sec:problem}) in which only first order deterministic effects are accounted for. The second is a diffusion scaling (see \S \ref{sec:diff}) under which a refinement of the fluid solution is considered to account for stochastic second-order terms.

\subsection{Main results}\label{sec23}
Our paper focuses on solving the LPOP asymptotically under the large
population limiting regime. Our first-order analysis uses a
fluid-scaling and captures a the deterministic elements of the system,
not accounting for stochasticity. This type of analysis allows us to
identify a simple near-optimal scheduling rule and gives us useful
insights about the original finite-population problem. Our
second-order analysis incorporates the stochastic elements back into
the model and offers a solution that is asymptotically optimal under
diffusion scaling for large time horizon $H$. In both scaling regimes we also identify the SG defined as the appropriately scaled difference between our proposed solution and the solution under complete information.

\subsubsection{Fluid scale}
Assume that the cost parameters are of the form postulated in
(\ref{eq:scaled_c}), and recall that this form was selected in such a
way that the two additive components of $J_n(\cdot)$ in
(\ref{eq:FPOP-n}) are both of order ${\mathcal
  O}(1)$. In the fluid scaling, the stochastic
  variation is `washed out,' ensuring that the stochastic optimization
  in (14) is approximated by a variational problem, in the large
  population asymptotic. This allows us to focus on the $\mathcal
  O(1)$ terms in the optimization.

Consider
\begin{equation}
E^f_n(t) =
\left\{
\begin{array}{ll}
1+\left\lfloor\frac{n\mu t}{p}\right\rfloor & t<H, \\
[0.2cm]
N_n & t=H,
\end{array}
\right.
\end{equation}
and its corresponding schedule
\begin{equation}
T^f_{i,n}= \min\left\{\frac{p}{n\mu}(i-1),H\right\},~~i=1,...,N_n,~~n\in {\mathbb N}.
\end{equation}
Then, we establish that $\{T^f_{i,n}\}$ is asymptotically optimal in the fluid scale. Specifically, we show that
\[
\lim_{n\rightarrow\infty} J_n\left(\left\{T^f_{i,n}\right\}\right)=\lim_{n\rightarrow\infty} V_n\doteq \bar{V}.
\]

The schedule $\{T^f_{i,n}\}$ satisfies the following properties:
\begin{itemize}
\item The appointment times up to time $H$ are at intervals of equal duration of $\frac{p}{n\mu}$ time units.
\item All patients who do not get appointments before time $H$ are scheduled to arrive at that time.
\item The arrival rate during $[0,H)$ is equal to the service rate of $n\mu$. Thus, it is asymptotically optimal to operate the system at a critically loaded heavy-traffic regime. Note that heavy traffic is obtained here as a {\em result} of optimality and not as an {\em assumption}.
\item {The critically loaded regime implies that, in the fluid scale,
  the server idle-time is negligible compared to the overage time, when all customers have been served, and that no customers wait during $[0,H)$.}
\end{itemize}

In terms of SG, our results show that in the fluid scaling this gap vanishes in the limit. Specifically, we show that
\begin{equation}
\lim_{n\rightarrow\infty}\gamma_n=\lim_{n\rightarrow\infty}(V_n-V_n^{\rm CI})=0.
\end{equation}
This result implies that, at the fluid scaling, knowing whether customers will show or not and their actual service time is only marginally beneficial to the system manager. In particular, knowing these quantities on average is sufficient at the fluid level.

\subsubsection{Diffusion scale}
In diffusion scaling we are interested in fluctuations about the fluid solution that are of order ${\mathcal O}(1/\sqrt{n})$. Specifically, we focus on the centered and scaled cost function
\begin{equation}\label{43}
\hat{J}_n\left(\left\{T_{i,n}\right\}\right)=\sqrt{n}\left(J_n\left(\left\{T_{i,n}\right\}\right)-\bar{V}\right),
\end{equation}
and its corresponding centered and scaled value function
\[
\hat{V}_n = \inf_{\left\{T_{i,n}\right\}\in\calT_n} \hat{J}_n\left(\left\{T_{i,n}\right\}\right).
\]
Our diffusion scale results are for a large time horizon. To state these results we will add the time horizon $H$ as a subscript to all relevant quantities.
Consider the following scheduling function
\begin{equation}
E^d_{n,H}(t) =
\left\{
\begin{array}{ll}
1+\left\lfloor\frac{nt}{p}\left(\mu-\frac{c^*_H}{\sqrt{n}}\right)\right\rfloor & t<H, \\
[0.2cm]
N_n & t=H,
\end{array}
\right.
\end{equation}
where
$c^*_H=\sqrt{\frac{c_w\left(p(1-p)+\mu^3\sigma^2\right)}{2\left(c_w(\bar{\tau}-H)+c_o/\mu\right)}}$,
with $\bar{\tau}=\frac{p\alpha}{\mu}$,
and its corresponding schedule
\begin{equation}
T^d_{i,n,H}= \min\left\{\frac{p}{n\left(\mu-c^*_H/\sqrt{n}\right)}(i-1),H\right\},~~i=1,...,N_n,~~n\in {\mathbb N}.
\end{equation}
Then, we establish that $\{T^d_{i,n,H}\}$ is asymptotically optimal in the diffusion scale. More precisely, under the assumption that the service times $\nu_i$ possess a $3+\eps$ moment,
we show that
\[
\lim_{H\rightarrow\infty}\limsup_{n\rightarrow\infty} \frac{1}{H} \hat J_{n,H}\left(\left\{T^d_{i,n,H}\right\}\right)=\lim_{H\rightarrow\infty}\liminf_{n\rightarrow\infty} \frac{1}{H}\hat{V}_{n,H}.
\]

The schedule $\{T^d_{i,n,H}\}$ satisfies the following properties:
\begin{itemize}
\item The appointment times up to time $H$ are at intervals of equal duration of $\frac{p}{n\left(\mu-c^*_H/\sqrt{n}\right)}$ time units. This interval duration deviates from the fluid schedule by a term of the order of ${\mathcal{O}}(1/\sqrt{n})$
\item In the scaling limit, the queue length process converges to a reflected Brownian motion (RBM) on the half-line, with a constant negative drift $-c^*_H$.
\item The constant $c^*_H$ is obtained by minimizing a cost with two additive terms; one that is proportional to $c^*_H$, and represents that server idleness cost, plus  a term that is inversely proportional to $c^*_H$, that represents the holding cost.
\end{itemize}

In terms of SG, it turns out that while in the fluid scale the SG is negligible, it is strictly positive in the diffusion scale. Specifically, let $\hat{\gamma}_{n,H}=\sqrt{n}\left(V_{n,H}-V^{CI}_{n,H}\right)$. Then, we show that
\[
\lim_{H\rightarrow\infty}\liminf_{n\rightarrow\infty}\frac{1}{H}\hat\gamma_{n,H}=
\lim_{H\rightarrow\infty}\liminf_{n\rightarrow\infty}
\frac{1}{H}\hat{V}_{n,H} =: \hat V^*> 0,
\]
{where $\hat V^*$ is a constant.}

\subsection{Exact analysis of the CI problem}\label{sec24}

In the CI problem, all stochastic data are known to the decision maker.
It thus can be treated, for each realization of these data,
as a deterministic allocation problem.
The number of show ups is given by $\Xi(N)$, and the total amount
of work associated with them is $\tau^*:=V(\Xi(N))=\inf\{t:S(t)\ge\Xi(N)\}$.
Hence $\tau^*$ is a lower bound on $\tau$ for any CI-schedule.
It is easy to see that there is no gain by allowing the system to
be empty (thus the server idle) for some time
prior to the time of completion of all jobs, $\tau$.
Indeed, if there is any interval $[a,b)$ on which the system is empty, and there
is still a job (or more) that will show up at time $b$, then advancing all jobs
scheduled at $b$ or later by $b-a$ units of time (so that the job originally scheduled
at $b$ arrives at $a$, etc.) does not affect the waiting time of any of the jobs
arriving at times $t\ge a$, and it can only decrease the overage time.
Thus it suffices to consider only allocations for which $B(t)=t$,
(equivalently, $Q(t)\ge1$) for all $t<\tau$.
Moreover, for all such allocations, clearly $\tau=\tau^*$.
As a result, \eqref{eq:queue-length} gives
\begin{equation}\label{56}
Q=A-S \text{ on } [0,\tau^*].
\end{equation}
The sample path of $A(t)$, $t\in[0,\iy)$,
can be any member of $\bar D[0,\iy)$ that is integer valued and satisfies
\begin{equation}\label{54}
A(t)=\Xi(N) \text{ for all } t\ge H.
\end{equation}
Now, in view of \eqref{56}, the requirement $Q(t)\ge1$ alluded to above implies
\begin{equation}\label{55}
A(t)\ge1+S(t) \text{ for all } t<\tau^*.
\end{equation}
Among all paths satisfying \eqref{54} and \eqref{55}
there is one that is pointwise minimal, namely
\[
A^*(t)=\begin{cases} 1+S(t), & t<H,\\ \Xi(N), & t\ge H.
\end{cases}
\]
It corresponds to allocating the jobs in such a way that there are no
waiting customers during $[0,H)$, and, in the event that $\tau^*\ge H$,
the remaining customers, that are $\Xi(N)-(1+S(H-))$ in number,
are all scheduled at the very last moment, $H$.
By \eqref{eq:waiting-time} and \eqref{56}, $W=\int_0^{\tau^*}(A(s)-S(s)-1)^+ds$.
Thus $W$ is monotone in $A$ in the following sense: If $A(t)\ge \tilde A(t)$
for all $t\in[0,\tau^*]$ then $W\ge\tilde W^*$. We conclude that
$A^*$ minimizes $W$, and since we have already mentioned that it minimizes
$\tau$, it also minimizes their weighted sum, in pathwise sense.
Consequently, this CI-schedule minimizes the cost \eqref{52}.

Finally, we can also compute this cost. On the event $\tau^*<H$, $W^*=0$;
when $\tau^*\ge H$,
\[
W^*=\int_H^{\tau^*}(\Xi(N)-S(s)-1)^+ds=\int_H^{\tau^*}(\Xi(N)-S(s)-1)ds.
\]
We thus obtain
\begin{align}\notag
  \label{57}
  V^{\rm CI}
  &=c_w\E[W^*]+c_o\E[(\tau^*-H)^+]\\
  &=c_w\E\Big[\int_H^{H\vee\tau^*}(\Xi(N)-S(s)-1)ds\Big]
  +c_o\E[(\tau^*-H)^+].
\end{align}

\section{Large population asymptotics: fluid scale} \label{sec:problem}
This section studies our appointment scheduling model in fluid scale.
By exploiting the stochastic regularity
that emerges in this scaling limit, we identify a
deterministic, first-order approximation to the FPOP that
governs the limit behavior. We find an optimal solution to this limiting problem and show that the cost associated with this solution is asymptotically achievable
in the fluid scale limit.

We start by stating and solving a formal fluid problem. Later we show that the optimal value of this fluid problem constitutes a lower bound on the fluid-scaled FPOP. Subsequently, we show that this value also constitutes an upper bound on the fluid-scaled  FPOP shown by identifying a sequence of simple policies for the FPOP that asymptotically achieve this value.  Thus, we establish that this sequence of policies is asymptotically optimal in the fluid scale.

\subsection{Fluid Model} \label{sec:fluid-model}

Let $\{E_n\}$ be an arbitrary sequence
of scheduling functions.
Following~\eqref{eq:arrival}, the fluid-scaled
cumulative arrival process is defined as:
\[
\bar A_n = \frac{1}{n} \Xi \circ n \bar E_n,
\]
where $\bar E_n = \frac{1}{n} E_n$ is the fluid-scaled
schedule. The FLLN implies that, as $n\to\iy$,
\begin{equation} \label{eq:arrival-count-flln}
\frac{1}{n} \Xi(\lfloor n e \rfloor)
\To pe,
\end{equation}
where $e : \bbR_+ \to \bbR_+$ is the identity map; note that here the
convergence is u.o.c., even though the pre-limit processes are assumed
to exist in the Skorokhod $J_1$ topology, since the limit process is
continuous. This will be the case in the remainder of the discussion,
unless noted otherwise. Throughout this section we use the notation $\eps_n$ for a generic sequence of stochastic
processes that converge to the zero process in probability as $n \to
\infty$, as well as for a generic sequence of RVs that converges to zero in probability.
It follows from \eqref{eq:arrival-count-flln} that
\begin{equation}
\bar A_n = n^{-1} \Xi(n \bar E_n) = p \bar E_n + \eps_n.
\label{eq:arrival-flln}
\end{equation}
Also, the FLLN for renewal processes \cite[Chapter 5]{ChYa01}
implies that, as $n\to\iy$,
\begin{equation}
  \label{eq:service-flln}
  \bar S_n := \frac{1}{n} S_n \To \m e.
\end{equation}
As a result, $\bar S_n(B_n(t)) = \mu B_n(t) + \eps_n$.
In view of these identities, the fluid-scaled queue length process $\bar Q_n := n^{-1} Q_n$, is given by

\begin{eqnarray*}
  \bar Q_n &=& \bar A_n - \bar S_n \circ B_n\\
  &=& p \bar E_n - \mu e + \mu (e -  B_n) + \eps_n.
\end{eqnarray*}
The idleness process $I_n := e -  B_n$ is non-decreasing
and is flat on excursions of $\bar Q_n$ away from zero,
and thus $(\bar Q_n,\mu I_n)$ forms
a solution to the Skorohod problem
with data $p\bar E_n-\mu e+\eps_n$.
That is,
\begin{equation} \label{eq:queue-flln}
\left( \bar Q_n, \mu I_n\right) = \left( \Gamma_1(p \bar E_n - \mu e +
  \eps_n), \Gamma_2\left( p \bar E_n - \mu e + \eps_n\right) \right),
\end{equation}
with $\Gamma_1(x) = x + \Gamma_2(x)$ and $\Gamma_2(x)(t) = \sup_{0 \leq s
  \leq t} (- x(s))^+$.

Recall the fluid-scaled makespan RV
\(
\bar W_n := n^{-1} W_n = n^{-1}\int_0^\infty \left( Q_n(t) -1 \right)^+ dt,
\)
and the overage time $O_n = (\tau_n - H)^+$, where $\tau_n = \inf \{t
> 0 : S_n(B_n(t)) \geq
\Xi(N_n)\}$. From~\eqref{eq:arrival-flln},~\eqref{eq:service-flln} and
\eqref{eq:queue-flln},
by formally removing the error terms, we derive a fluid model
as follows.
Let $\calL=\{\la\in\sDplus[0,\iy):\la(t)=\al\text{ for } t\ge H\}$. Given $\la\in\calL$,
let $a=p\la$ and $q = \Gamma_1(a - \mu e) = a - \m e + \eta$, where
$\eta(t) := \Gamma_2(a - \m e)(t)$ is the correction term (or
Skorokhod regulator term). These are fluid models
for the arrival process and queue length, respectively.
Let also $\bar B(t) := t - \mu^{-1} \eta(t)$ stand
for cumulative busyness, and let $\tau = \inf \{t > 0 : \mu \bar B(t) \geq a(H)\}$,
$\bar W = \mu^{-1}\int_0^\infty q(t) da(t)$ and $\bar O =
(\tau-H)^+$ denote the fluid models for the termination time, wait and overage time,
respectively.
A \textit{fluid optimization problem} (FOP) is formulated by letting
\begin{equation}
  \label{eq:fluid-cost}
  \bar J(\l) = c_w \bar W + c_o  \bar O,
\end{equation}
and
\begin{equation}
  \label{eq:FOP}
  \bar V := \inf_{\lambda \in \calL} \bar J(\lambda).
\end{equation}

In \S \ref{sec:fluid-schedule} below we show that there exists
a $\l^*\in\calL$ that attains the minimum in \eqref{eq:FOP}.
In \S \ref{sec:fluid-lower-bound} we show
that the FOP value $\bar V$ is an asymptotic lower bound on the fluid scale cost
$ J_n(\{T_i\})$ under an arbitrary sequence of schedules.
In \S \ref{sec:fluid-upper-bound} we construct a bespoke sequence of
finite population schedules that asymptotically \textit{achieves}
$\bar V$, thus proving asymptotic optimality in fluid scale.

\subsection{Fluid Optimal Schedule} \label{sec:fluid-schedule}
Under our overload assumption \eqref{eq:overload}, it is straightforward to see that
our fluid model satisfies, for any $\la\in\calL$, $0<\tau - H = \mu^{-1}
q(H)$. An optimal control should minimize the tradeoff
between the fluid overage time and fluid makespan.
The main result of this section in Proposition \ref{thm:fop} identifies such
an optimal control.

\begin{proposition} \label{thm:fop}
A fluid optimal control $\lambda^* \in \calL$ is given by
  \begin{equation}
    \label{eq:fop}
    \lambda^*(t) = \begin{cases}
      p^{-1} \m t ~&~ t \in [0,H)\\
      \a ~&~ t \geq H.
      \end{cases}
  \end{equation}
\end{proposition}
The proof is fairly straightforward and is a consequence of the
following lemma that, loosely speaking, argues that it is optimal for the server to be busy at all times until all work has been completed.

\begin{lemma}\label{lem:non-idle}
For any control $\l \in \sL$ with the corresponding correction term $\eta(t) \geq 0$, there exists a control
$\tilde{\lambda} \in \sL$ such that the corresponding $\tilde \eta(t) = 0$ for all $t
\in [0,\tilde \tau]$, where $\tilde \tau = \inf \{t > 0 : \m t \geq p
\tilde \lambda(H) \}$, and $\bar J(\l) \geq \bar J (\tilde \lambda)$.
\end{lemma}

This result implies that
the cost associated with any control where the correction term is non-zero in $[0,H]$ can be
improved by choosing a control where the correction term is
zero. Thus, if an optimal control exists, there exists one such control $\l$ that satisfies $p (\l(t_2) - \l(t_1)) \geq \m
(t_2 - t_1)$, where $0 \leq t_1 \leq t_2 \leq H$. We now prove Proposition~\ref{thm:fop}.

\skp

\noi{\bf Proof of Proposition~\ref{thm:fop}.}
Following Lemma~\ref{lem:non-idle}, we will assume that $\l \in \sL$ is
such that the correction term $\eta$ is zero in $[0,H]$. Let $\bar c_w = p c_w (\a \mu)^{-1}$ and $\bar c_o = c_o (\mu)^{-1}$. Integrating the cost function $\bar J(\lambda)$ by parts and using the fact that $\lambda(H)=1$ and $q(0-) = 0$,
\begin{align*}
\bar J(\la)&=\bar c_w\Big[q(t)\la(t)|_0^H-\int_0^H\la(t)dq(t)\Big]+\bar c_oq(H)\\
&=\bar c_w\Big[q(H)-\int_0^H\la(t)(pd\la(t)-\mu dt)\Big]+\bar c_oq(H)\\
&= -\bar c_w\frac{\a^2 p}{2}+\bar c_w\mu\int_0^H\la(t)dt
    +(\bar c_w+\bar c_o)(p \a-\mu H).
\end{align*}
The only term that can be controlled is the second one, which is optimized by
\[
\la(t)=\begin{cases}
  p^{-1} \mu t, & t\in[0,H),\\
  \a, & t\ge H,
\end{cases}
\]
since $p (\l(t_2) - \l(t_1)) \geq \m(t_2 - t_1)$ by the discussion following Lemma~\ref{lem:non-idle}.
\qed

\skp

\noi{\bf Proof of Lemma~\ref{lem:non-idle}.}
Recall that fluid cost under $\l$ is
\[
\bar J(\la) = c_w p\mu^{-1}\int_0^H q(t) d \la (t) + c_o (\tau - H),
\]
where we have made use of fact that $a = p \la$ and $d a(t) = 0$ for
all $t > H$. Now, let
$\tilde \l := \l + p^{-1} \eta$ on $[0,H)$ and $\tilde \l(H) := \a$, so that $\tilde \l \in \sL$. Observe that the queue length under
$\tilde \l$ is
\[
\tilde q(t) =
\begin{cases}
  q(t) & \text{if}~t \in [0,H)\\
  \a - \m t & \text{if}~t \in [H,\tilde \tau],\\
  0 & ~\forall~t > \tilde \tau.
\end{cases}
\]
 Notice that $\tilde
q(H) \leq q(H) = \a - \m H + \eta(H)$, implying that $\tilde \tau \leq
\tau$. The associated fluid cost is
\[
\bar J(\tilde \l) = c_w p \mu^{-1}\int_0^H q(t) d\tilde \l(t) + c_o (\tilde \tau-H).
\]
Since $ q(t) d \eta(t) = 0$, it follows that
\[
\bar J(\tilde \l) = c_w p \mu^{-1}\int_0^H q(t) d \l(t) + c_o (\tilde \tau - H) \leq \bar J(\l).
\]
\qed

\skp

Some remarks on Proposition~\ref{thm:fop} are warranted. First, and most
importantly,
under~\eqref{eq:overload} the optimal
schedule ensures that the queue length is zero
in $[0,H)$ and positive in $[H,\tau]$ and zero for $t \geq \tau$. This is an intuitively satisfying result
in the sense that, given the system operator's goal of minimizing the
makespan while simultaneously minimizing the overage time, it would
make most sense to fully utilize the available capacity but
\textit{not} overload the system. Thus, the optimal schedule matches the
 arrival `rate' with the effective service rate $p^{-1} \mu$ in the
 interval $[0,H)$, and schedules the remainder (of $\a - p^{-1}\mu H$ fluid
 units) at $H$. In other words, the heavy-traffic
condition emerges as a \textit{consequence} of optimization. This
result is in stark contrast to most queueing control problems where
the heavy-traffic condition is assumed at the outset.
~Notice too that the fluid optimal solution parallels the sample
path-wise solution in the complete information problem. In the absence
of stochastic variation in the traffic and service, it is clearly
possible to optimally arrange the traffic such that there is no waiting
for any of the jobs who do turn up, and schedule the remainder at the
end of the horizon.

\subsection{Lower bound on the fluid scale cost}\label{sec:fluid-lower-bound}
Recall that the fluid scaled makespan $\bar W_n$ is defined as
\(
\bar W_n = n^{-1}\int_0^\infty (Q_n(s) -1)^+ ds.
\)
Now, fix $0<K<\infty$ (to be determined later) and note that
\begin{eqnarray}
\nonumber
\bar W_n &\geq& n^{-1}\int_0^K (Q_n(s)-1)^+ ds\\
\label{eq:makespan-approx}
&\geq& \int_0^K \bar Q_n(s) ds - K n^{-1}.
\end{eqnarray}
By \eqref{eq:queue-flln}, there
exists a sequence of processes $\{\eps_n\}$ that converges to the zero
process in probability as $n \to \infty$, such that
\[
\bar Q_n  = \Gamma_1(p \bar E_n - \m e + \eps_n).
\]
By the Lipschitz continuity of the Skorokhod regulator map it follows that
\[
\bar Q_n \geq \Gamma_1(p \bar E_n - M) + \eps_n,
\]
where $M=\mu e$.
Substituting this into~\eqref{eq:makespan-approx} we observe that
\begin{equation}
  \nonumber
  \bar W_n \geq \int_0^K \Gamma_1(p \bar E_n-M)(s) ds + \eps_n - K n^{-1},
\end{equation}
clearly implying that
\begin{equation} \label{eq:makespan-bound}
\bar W_n \geq \inf_{\l \in \sDplus[0,\infty)} \int_0^K \Gamma_1(p \l - M)(t) dt + \eps_n - K n^{-1}.
\end{equation}

From the definition of $\tau_n$~\eqref{50-n} it follows that
\[
\tau_n \geq \inf\{t : \m B_n(t) + \eps_n \geq \a p\},
\]
where we have used the fact that $\bar S_n(B_n(t)) = \m B_n(t) +
\eps_n$ and $\bar A_n(H) = \a p + \eps^{'}_n$; note that we use $\eps_n$ to
represent the difference between the two mean-zero error sequences. Since $t \geq B_n(t)$
we have
\(
\tau_n \geq \inf\{t : \mu t + \eps_n \geq \a p \},
\)
implying that $\tau_n \geq  \a p \m^{-1} + \eps_n = \bar{\tau}  + \eps_n$. It follows that
\begin{equation}
  \label{eq:overage-bound}
  (\tau_n - H)^+ \geq (\bar{\tau} - H + \eps_n)^+.
\end{equation}

Let $\bar{\tau} = p \a \mu^{-1}$ and consider
\[
\tilde V := c_w
\inf_{\la\in \sL} \left\{ \int_0^\iy
\Gamma_1(p \lambda - M)(t) dt \right\} + c_o (\bar{\tau} - H)^+.
\]
As the
next lemma shows $\tilde V$ equals the FOP value, and can
be achieved by the optimal schedule in Proposition~\ref{thm:fop}.
\begin{lemma}\label{lem:optimal-fluid-properties}
  \noi 1.
\(
\bar V=\tilde V.
\)
\\
\noi
2. The upper limit $\iy$ in the integral can be replaced by any $K$ sufficiently large.

\noi
3. The minimum is attained by $\l^*$ defined in~\eqref{eq:fop}.
\end{lemma}

\noi{\bf Proof of Lemma ~\ref{lem:optimal-fluid-properties}.}
Consider $\tilde V$ first and recall that for a fixed $\l$, $\tau
= \inf \{t > 0 : \mu (t - I(t)) = \a p\}$. If $\tau > \bar{\tau}$, it automatically follows that $I(\bar{\tau}) > 0$.  This implies that the makespan cost $c_w\int_0^\infty (p \l(t) - \mu (t - I(t))) dt$ is not optimal. To see this, note that the makespan cost can be lower bounded by choosing $\l'$ such that $\l'(t) = \mu t$ for all $t \in [0,\tau]$, and $\tau = \bar{\tau}$ in this case. Thus, any optimal solution should be such that $I(t) = 0$ up to $\tau=\bar{\tau}$, in which case we minimize
\[
\int_0^{\bar{\tau}} q(t) dt = \int_0^{\bar{\tau}} p \l(t) dt - \int_0^{\bar{\tau}} \mu t dt,
\]
where only the first term on the right hand side is controlled. The minimizing schedule $\l$ that satisfies the constraint that $\l(t) = \a,~t\geq H$ is
\[
\l(t) =
\begin{cases}
  p^{-1} \mu t & t \in [0,H)\\
  \a & t \geq H.
\end{cases}
\]
\qed

\skp

Now we will show that Lemma~\ref{lem:optimal-fluid-properties} together
with~\eqref{eq:makespan-bound} and~\eqref{eq:overage-bound} implies
that the FOP value lower bounds the fluid-scaled cost.

\begin{theorem}\label{prop:fluid-cost-lower-bound}
  The fluid-scaled cost of an arbitrary sequence of schedules,
  $\{T_i\}$, is asymptotically lower bounded by the fluid optimal
  value $\bar V$. That is,
  \[
  \liminf_{n \to \infty} \bar J_n(\{T_i\}) \geq \bar V.
  \]
\end{theorem}

\noi{\bf Proof of Theorem ~\ref{prop:fluid-cost-lower-bound}.}
Let $j_n := c_w \bar W_n + c_o O_n$ represent the random
cost incurred by following schedule
$\{T_i\}$. Since the constant $K$ was arbitrary, we can set it to be
greater than $\bar{\tau}$. Lemma ~\ref{lem:optimal-fluid-properties} together
with~\eqref{eq:makespan-bound}
and~\eqref{eq:overage-bound} imply that,
\begin{equation*}
j_n \geq (\bar V + \eps_n - K n^{-1}) \vee 0.
\end{equation*}
Observe that $j_n \To \bar V$ as $n \to \infty$. Since $j_n \geq 0$,
Fatou's Lemma~\cite{EtKu2008} implies
that
\[
\liminf_{n \to \infty} \bbE[j_n] \geq \bbE[\liminf_{n\to\infty} j_n] =
\bar V.
\]
\qed

\skp

\subsection{Upper bound on the fluid scale cost}\label{sec:fluid-upper-bound}
We now construct a sequence of scheduling policy whose fluid-scaled cost is
asymptotically upper-bounded by the FOP value $\bar V$. Given the lower bound result in the previous subsection, this sequence is thus asymptotically optimal in the fluid limit.

Recall, that the fluid optimal schedule is
\[
\lambda^*(t) = \begin{cases}
\mu p^{-1}t & \forall t \in [0,H)\\
\a & \forall t \geq H.
\end{cases}
\]
Consider the following sequence of scheduling functions indexed by $n$
\begin{equation*}
  E_n^f(t) := \begin{cases}
1 + \left\lfloor\frac{n \m t}{p} \right\rfloor & t < H,\vspace{0.2cm}\\
N_n & t=H,
\end{cases}
\end{equation*}
and its corresponding schedule
\begin{equation}
T^f_{i,n}= \min\left\{\frac{p}{n\mu}(i-1),H\right\},~~i=1,...,N_n,~~n\in {\mathbb N}.
\end{equation}
The scheduling function $E_n^f$ is interpreted as follows: for each $n$
customers are scheduled to arrive one-at-a-time at
uniformly spaced intervals of length $p(n\mu)^{-1}$ up to time $H$, with the leftover $N_n -(1+\left\lfloor n\mu H/p\right\rfloor)$ customers who are scheduled to arrive at time $H$.

The main result of this section establishes the fact that
the expected fluid-scaled cost $\bar J_n (\{T_{i,n}^f\})$ converges to the fluid optimal value as well.

\begin{theorem}~\label{thm:asymptotic-achievability}
  For each $n$, suppose that traffic schedule is $\{T^f_{i,n}\}$. Then
  \begin{equation}
    \limsup_{n \to \infty} \bar J_n(\{T_{i,n}^f\}) \leq \bar V.
  \end{equation}
\end{theorem}

\begin{lemma}~\label{lem:finite-pop-sched}
  The finite population schedule satisfies
  \(
    \E^f_n \to \l^* ~\text{uniformly on compacts as}~n\to\infty.
    \)
\end{lemma}

 Now, let $(Q_n,I_n)$ represent the queue length and the idleness
processes when traffic is scheduled per $E^f_n$. FLLN's for the arrival
and service processes, and Lemma~\ref{lem:finite-pop-sched},
together imply the following result:

\begin{lemma}\label{lem:upper-bound-convergence}

\noindent (i) The fluid-scaled queue length and idleness processes
satisfy an FLLN: $(\bar Q_n, \bar I_n) \To (q^*, \iota^*)$ as $n
\to \infty$ where
\[
q^*(t) =
\begin{cases}
  0 & t \in [0,H)\\
  \a p - \mu t & t \in [H,\bar{\tau}]\\
  0 & t > \bar{\tau},
\end{cases}
\]
and
\[
 \iota^*(t) =
\begin{cases}
  0 & t \in [0,\bar{\tau}]\\
  t-\bar{\tau} & t > \bar{\tau}.
\end{cases}
\]

\noindent (ii) The fluid-scaled makespan and overage time RVs
satisfy: $\left( \bar W_n, O_n \right) \To \left(
  \bar W^*, \bar O^* \right)$ as $n \to \infty$, where $\bar W^* = p\mu^{-1}
\int_0^\infty q^*(t) d\l^*(t)$ and $\bar O^* = (\bar{\tau} - H)$.
\end{lemma}

\noi{\bf Proof of Theorem \ref{thm:asymptotic-achievability}.}
Let $J^R_n(\{T^f_{i,n}\}) = c_w \bar W_n + c_o  O_n$ be the random cost. We start by noting that the convergence result
in~Lemma~\ref{lem:upper-bound-convergence} implies that the random cost $J^R_n(\{T^f_{i,n}\}) \to \bar J(\l^*)$ weakly converges to
 $\bar V$ as $n \to
\infty$. \nocite{EtKu2008} implies that $ J_n(\{T^f_{i,n}\}) =
\bbE[J^R_n(\{T^f_{i,n}\})]$ will converge to $\bar V$ provided that $J^R_n(\{T^f_{i,n}\})$ is
uniformly integrable. The remainder of this proof is dedicated to
proving this claim.

We prove that $J^R_n(\{T_{i,n}^f\})$ is uniformly integrable by showing that $\bbE| J^R_n(\{T_{i,n}^f\})|^2
\leq C < \infty$ for all $n \in \bbN$. Consider the
sequence $\{O_n\}$ first. Note that it suffices to consider the case
where $\tau_n > H$. The number of jobs waiting in the queue at the end
of the arrival horizon $H$ is $N_n-D_n(H) > 0$, where $D_n(H)$ is the
number of departures in $[0,H]$. Since there are no
more arrivals after time $H$ it can be seen that
\begin{equation} \label{eq:bound1}
\tau_n - H \leq \sum_{i=D_n(H)+1}^{N_n} \nu_i^n \leq
\sum_{i=1}^{N_n} \nu_i^n.
\end{equation}

Now, let $\Upsilon(m) := \sum_{i=1}^m \nu_i$ and $\Upsilon_n(m) := \sum_{i=1}^m
\nu_i^n$. Minkowski's inequality implies that $\left(\bbE
\left|\Upsilon(N_n)\right|^2 \right)^{1/2}
\leq N_n \left( \bbE \left| \Upsilon(1) \right|^2 \right)^{1/2}$. Therefore, we obtain
\begin{equation} \label{eq:workload-bound}
\left( E \left|\Upsilon_n(N_n)\right|^2 \right)^{1/2} \leq \frac{N_n}{n}
\left( E \left|\Upsilon(1)\right|^2 \right)^{1/2}
\leq \a \left( \bbE \left|\Upsilon(1)\right|^2 \right)^{1/2},
\end{equation}
where the last inequality follows from the fact that $N_n/n \leq \a$
by definition. Equation~\eqref{eq:bound1} and this bound imply
that
\[
\bbE|(\tau_n - H)^+|^2 \leq \bbE|\tau_n - H|^2
\leq \a \bbE \left|V(1)\right|^2 < \infty,
\]
where the finiteness of the second moment is by assumption. Since the
bound is independent of $n$, it follows that $O_n=(\tau_n - H)^+$ are uniformly integrable.

Now consider the fluid-scaled makespan. Using the fact that the queue
drains out and remains empty after $\tau_n$ it follows that
\[
\bar W_n = n^{-1} \int_0^H Q_n(t) dt + n^{-1} \int_H^{H \wedge
  \tau_n} Q_n(t) dt.
\]
Note that the first term on the right hand side of the inequality is
bounded above by $n^{-1}N_n H \leq \a H$. Thus, it suffices to
consider the second term when $\tau_n > H$. As there are $N_n - D_n(H)$
jobs waiting for service at the end of the horizon, it follows that
\begin{align*}
n^{-1}\int_H^{\tau_n} Q_n(t) dt &\leq n^{-1} \left\{ (N_n - D_n(H)) \nu_{D_n(H)+1}^n +
(N_n - D_n(H)-1) \nu_{D_n(H)+2}^n + \cdots + \nu_{N_n}^n\right\}\\
&= \frac{1}{n} \sum_{i=1}^{N_n-D_n(H)} (N_n - D_n(H) + 1 - i)
    \nu_{D_n(H)+i}^n\\
&= \frac{N_n -D_n(H)}{n} \sum_{i=1}^{N_n - D_n(H)} \nu_{D_n(H)+i}^n -
    \frac{1}{n} \sum_{i=1}^{N_n - D_n(H)} (i-1) \nu_{D_n(H)+i}^n\\
&\leq \frac{N_n}{n} \sum_{i=1}^{N_n} \nu_i^n,
\end{align*}
where the last inequality follows from the fact that $n^{-1}
\sum_{i=1}^{N_n-D_n(H)}(i-1)\nu_{D_n(H)+i}^n \geq 0$ and $D_n(H) \geq 0$ for all $n \geq 1$. Using the bound in~\eqref{eq:workload-bound} we have
\[
n^{-1}\left( \bbE \left (\int_H^{\tau_n} Q_n(t) dt \right)^2\right)^{1/2}
\leq \a \left( \bbE|\Upsilon(1)|^2 \right)^{1/2}.
\]
Thus, it follows that $\bbE|\bar W_n|^2$ is uniformly bounded for all $n
\in \bbN$, implying that $\{\bar W_n\}$ is uniformly integrable. Finally, since
$\{J^R_n(\{T_{i,n}^f\}),~n\geq 1\}$ is a sequence of RVs that are each linear
combinations of uniformly integrable RVs, it is uniformly integrable as well.
\qed

\skp

Thus, Theorem~\ref{thm:asymptotic-achievability} shows that the family of finite
population schedules $\{E^f_n\}$ is asymptotically optimal in the sense
that the FOP value can be achieved in the large population limit.

\skp

\noi{\bf Proof of Lemma \ref{lem:finite-pop-sched}.}
Fix $t \in [0,H)$ and $n \geq 1$. By definition it follows that
\[
  \Big| \frac{E_n^f(t)}{n} - \l^*(t) \Big| \leq \frac{1}{n} +
                                                       \frac{1}{n}
                                                       \Big(\frac{n
                                                       \m t}{p} -
                                                       \left\lfloor \frac{n
                                                       \m t}{p}\right\rfloor
                                                       \Big)
  \leq \frac{2}{n}.
\]
On the other hand, fix $t \geq H$ and observe that
\begin{equation}
  \Big| \frac{E_n^f(t)}{n} - \a \Big| \leq \Big| \frac{N_n}{n} -
    \a \Big| \leq \frac{2}{n}.
\end{equation}
These two bounds are independent of $t$, proving the lemma.
\qed

\skp

\noi{\bf Proof of Lemma \ref{lem:upper-bound-convergence}.}
Part (i) follows by using the FLLN's alluded to
in~\eqref{eq:arrival-count-flln} and~\eqref{eq:service-flln}, and Lemma~\ref{lem:finite-pop-sched}, and applying the continuous
mapping theorem to the tuple $(\bar Q_n, \bar I_n)$
(see~\eqref{eq:queue-flln}).

For part (ii), observe that by definition the fluid-scaled makespan also satisfies
\[
\bar W_n = n^{-1}\int_0^\infty W_n(t) d A_n(t),
\]
where $W_n(t) =
n^{-1}\sum_{i=1}^{\lfloor nt \rfloor} \nu_i - t + B_n(t)$ is the
workload process. Let $w^* := \mu^{-1} q^*$, then,
\begin{align*}
\nonumber
  &\left| \int_0^\infty W_n(t) d \frac{A_n(t)}{n} -
  \int_0^\infty w^*(t) da(t)  \right|\\
   &\quad\leq \left|
                                                          \int_0^\infty
                                                          W_n(t)
                                                          d\left(\frac{A_n(t)}{n}
                                                          -
                                                          a(t)\right)\right|
  + \left|\int_0^\infty \left( W_n(t) - w^*(t)\right) da(t) \right|\\
\nonumber
&\quad\leq \int_0^\infty W_n(t) |\gamma_n(dt)| + \int_0^\infty |W_n(t) -
  w^*(t)| da(t),
\end{align*}
where $\gamma_n := n^{-1}A_n - a$ is a signed measure.
Consider the second integral above and
observe that
\begin{align*}
  \int_0^\infty |W_n(t) - w^*(t)| da(t) \leq \sup_{0 \leq t \leq
  \tau_n \vee \bar{\tau}} |W_n(t) - w^*(t)|,
\end{align*}
where $\tau_n \Rightarrow \bar{\tau}$ as $n \to \infty$. Now, \cite[Proposition 3]{HoJaWa2015} implies
that $W_n \Rightarrow w^*$. Consequently, it follows that
\begin{align*}
  \int_0^\infty |W_n(t) - w^*(t)| da(t) \Rightarrow 0 ~\text{as}~n\to\infty.
\end{align*}

On the other hand, by the Jordan decomposition theorem~\cite[pp 346]{KoFo1970}
it follows that there exist positive, finite, measures $\gamma_n^+$ and
$\gamma_n^-$ such that $\gamma_n= \gamma_n^+ - \gamma_n^-$,
$\gamma_n^+([0,\infty)) < \infty$ and $\gamma_n^-([0,\infty)) < \infty$. Furthermore, the
(total) variation of the signed measure is $|\gamma_n|=\gamma_n^+ +
\gamma_n^-$. From the FLLN alluded to in~\eqref{eq:arrival-count-flln} it
follows that $|\gamma_n| \Rightarrow 0$ as $n \to \infty$. Then, by
\cite[Corollary 8.4.8]{Bo2007} it follows that
\begin{align}
\label{eq:decomp-converge}
\gamma_n^+ \Rightarrow 0~\text{and}~\gamma_n^- \Rightarrow 0 ~\text{as}~ n \to \infty.
\end{align}

Now, let $K \in (0,\infty)$ be a large
positive constant, so that $W_n^K(t) := n^{-1}\sum_{i=1}^{\lfloor nt
  \rfloor} \nu_i \wedge K - t + B_n(t) \leq n^{-1}\sum_{i=1}^{\lfloor nt
  \rfloor} \nu_i \wedge K \leq K$. Consider the integral,
\begin{align*}
  \int_0^\infty W_n^K |\gamma_n(dt)| &\leq K
                                       \int_0^\infty|\gamma_n(dt)|\\
  &= K \int_0^\infty(\gamma_n^+(dt) + \gamma_n^-(dt)) \Rightarrow 0
    ~\text{as}~n\to \infty,
\end{align*}
where the convergence follows from~\eqref{eq:decomp-converge}. Note
too that this convergence is true for any $0 < K < \infty$, implying
that
\[
\int_0^\infty W_n |\gamma_n(dt)| \Rightarrow 0 ~\text{as}~n \to \infty.
\]
\qed

\subsection{The stochasticity gap at the fluid scale}
The first-order deterministic FOP is solved by scheduling traffic to
match the available capacity. In addition, the previous two
sections have shown that the fluid scale cost is bounded by the FOP
value and proposed an asymptotically optimal schedule for the LPOP. At the same time, for the complete information (CI) problem, we identified, in \S \ref{sec24}, that a
$\Sigma_n$-measurable schedule~\eqref{57} that optimizes this
problem; this schedule allocates appointments such that there
are no waiting customers during $[0,H)$ and the server is never idle. Clearly, for the LPOP, there is a cost to be paid for
scheduling traffic without \textit{a priori} knowledge of the
randomness. The parallels between the CI and FOP optimal schedule and
Theorem~\ref{prop:fluid-cost-lower-bound} suggest that there may
be a gap between the LPOP value ($V_n$) and the value of the CI
problem ($V_n^{CI}$). We
quantify this \textit{stochasticity gap} by showing that $\gamma_n := V_n
- V_n^{CI} \geq 0$ decreases to zero as $n \to \infty$.
Recall that $\bar V = c_w \int_H^{\bar{\tau}} x(t) dt +
c_o (\bar{\tau} - H)$ is the value of the FOP. The following is the main result of this section
\begin{proposition} \label{thm:sg-fluid}
  The SG in the fluid limit is zero. That is,
  \(
  \gamma_n \to 0
  \)
as $n \to \infty$.
\end{proposition}

The proof of Proposition~\ref{thm:sg-fluid}
follows as a consequence of the following lemmas.

\begin{lemma}~\label{lem:overage-CI}
  The optimal overage time in the CI problem satisfies:
  \begin{equation}
    \label{eq:overage-CI}
    \bbE[(\tau_n^*-H)^+] \to (\bar{\tau} - H)~\text{as}~n \to \infty,
  \end{equation}
  where $\bar{\tau} = p (\a\mu)^{-1}$.
\end{lemma}

Now, let $X_n(t) := n^{-1}(\Xi(N_n) - S_n(t) - 1)^+$ and $x(t) = (\a p - \mu t)$ for all $t \geq 0$. We prove that the expectation of the integral $\int_H^{H\vee \tau_n^*} X_n(t) dt$ converges as $n \to \infty$.

\begin{lemma}\label{lem:makespan-CI}
  The optimal expected makespan of the CI problem satisfies
  \begin{equation}
    \label{eq:2}
    \bbE \left[ \int_H^{H\vee \tau_n^*} X_n(t) dt\right] \to \int_H^{\bar{\tau}} x(t) dt ~\text{as}~n \to \infty.
  \end{equation}
\end{lemma}

\skp

\noi{\bf Proof of Proposition \ref{thm:sg-fluid}.}
Note that $\gamma_n = V_n - V_n^{CI} = (V_n
- \bar V) + (\bar V - V_n^{CI})$. Then, Lemma~\ref{lem:overage-CI} and
Lemma~\ref{lem:makespan-CI} imply that $\bar V - V_n^{CI} \to 0$ as $n
\to \infty$. Theorem~\ref{thm:asymptotic-achievability} implies that $\limsup_{n
  \to \infty} (V_n - \bar V) \leq \limsup_{n\to \infty} (\bar J_n(\{T_{i,n}^f\}) -
\bar V) \leq 0$. On the other hand,
Theorem~\ref{prop:fluid-cost-lower-bound} implies that
$\liminf_{n\to \infty} (V_n - \bar V) \geq 0$. Thus, $ (V_n - \bar V)
\to 0$ as $n \to \infty$.
\qed

\skp

\noi{\bf Proof of Lemma ~\ref{lem:overage-CI}.}
By definition, $\tau_n^* := \Upsilon_n (\Xi(N_n))$, where recall that $\Upsilon_n(m) := \sum_{i=1}^m
\nu_i^n$. It is
straightforward to deduce that $(\tau_n^*-H)^+ \To  (\bar{\tau} - H)$ as $n \to
\infty$ as a consequence of the LLN. Since $\Xi(N_n) \leq N_n$ and $\Upsilon_n(\cdot) \geq 0$ for all $n \geq 1$, it follows that $(\tau_n^* - H)^+ \leq (\Upsilon_n(N_n) -
H)^+$. On the other hand, following~\eqref{eq:workload-bound} we have
\(
\bbE\left| (\Upsilon_n(N_n) - H)^+ \right|^2 \leq \a^2 \bbE\left| \Upsilon(1)^2,
\right| + H^2 < \infty,
\)
where recall that $\Upsilon(m) := \sum_{i=1}^m
\nu_i$,
implying that $\bbE\left| (\tau_n^* - H)^+ \right|^2
< \infty$. Therefore $(\tau_n^* - H)^+$ is uniformly integrable,
implying~\eqref{eq:overage-CI}.
\qed

\skp

\noi{\bf Proof of Lemma ~\ref{lem:makespan-CI}.}
We first prove that the optimal makespan converges in probability to
the limit on the right hand side. It suffices to assume that $\tau_n^* > H$, since we know from the proof of  Lemma~\ref{lem:overage-CI} that $\tau_n^* \To \bar{\tau}$ as $n \to \infty$ and that $\bar{\tau} > H$ (by the overload assumption \eqref{eq:overload}). Thus, consider
\begin{equation}\label{42}
\left | \int_H^{\tau_n^*} X_n(t) dt - \int_H^{\bar{\tau}} x(t) dt \right| \leq \left| \int_H^{\tau_n^*} X_n(t) dt - \int_H^{\bar{\tau}} X_n(t) dt \right| + \left| \int_H^{\bar{\tau}} X_n(t) dt - \int_H^{\bar{\tau}} x(t) dt\right|.
\end{equation}
Consider the first term on the RHS above, and observe that
\begin{align*}
\left| \int_H^{\tau_n^*} X_n(t) dt - \int_H^{\bar{\tau}} X_n(t) dt
  \right|&=\left| \int_{\tau_n^*}^{\bar{\tau}} X_n(t) dt \mathbf
  1_{\{\tau_n^* \leq \bar{\tau}\}} + \int_{\bar{\tau}}^{\tau_n^*} X_n(t) dt \mathbf
  1_{\{\tau_n^* > \bar{\tau}\}}\right| \\
  &\leq \alpha \left|\tau_n^* - \bar{\tau} \right|,
\end{align*}
where the last inequality follows from the fact that $X_n(t) \leq
\alpha$ for all $t \in [0,\infty)$. Therefore, it follows that $\left|
  \int_H^{\tau_n^*} X_n(t) dt - \int_H^{\bar{\tau}} X_n(t) dt\right| \To
0$ as $n \to \infty$. Next, consider the second term on the RHS. Using the facts that
$X_n(t) \To x(t)$ as $n \to \infty$ pointwise and $|X_n(t) -
x(t)| \leq \alpha$ for all $t \in [0,\infty)$, the bounded convergence
theorem implies that $\int_H^{\bar{\tau}} |X_n(t) - x(t)| dt \To 0$ as $n
\to \infty$. Thus, it follows that $\int_H^{H\vee \tau_n^*} X_n(t) dt
\To \int_H^{\bar{\tau}} x(t) dt$ as $n \to \infty$. Finally, observe that
$\left|\int_H^{H\vee \tau_n^*} X_n(t) dt \right|^2 \leq \a^2
|(\tau_n^* - H)^+|^2$, and from the analysis in
Lemma~\ref{lem:overage-CI} it follows that $\bbE \left|\int_H^{H\vee
    \tau_n^*} X_n(t) dt \right|^2 < \infty$. Therefore, the sequence
of integrals are uniformly integrable, implying~\eqref{eq:2}.
\qed

\section{Large population optima: diffusion scale} \label{sec:diff}

Some of the results in this section require a strengthening of the second moment
condition of the service time.
\begin{assumption}
  \label{assn1}
  The service times $\nu_i$ possess a finite $3+\eps$ moment, that is,
  $E[\nu_1^{3+\eps}]<\iy$ for some $\eps>0$.
\end{assumption}

\subsection{Model equations and BOP derivation}\label{sec41}

As captured by the fact that the FOP
is deterministic, the inherent stochasticity in the FPOP degenerates in fluid scale.
A more realistic setting should capture the effect
of the stochastic variation introduced by the no-shows and the random service times.
In this section we consider the scheduling
problem at diffusion scale that incorporates second-order effects.

Our first goal in this subsection is to write down equations for the various quantities
of interest related to the diffusion scale problem, and then use these equations
to propose a Brownian optimization problem (BOP).

By the fluid scale analysis, the rescaled asymptotically optimal schedule
$E^f_n(t)$ converges to $\la^*(t)=(p^{-1}\mu t)\wedge\alpha$.
Under our assumption~\eqref{eq:overload}, $p^{-1}\mu H < \alpha$, and so the function $\la^*$ has a jump of size $\la^H:=\alpha-p^{-1}\mu H$ at
$H$.
In general, we shall denote with a superscript $H$
quantities that correspond to the allocation at the singular time point $H$.

Denote
\[
\hat E_n(t)=\frac{E_n(t)-np^{-1}\mu t}{\sqrt n},\qquad t\in[0,H).
\]
Recall that $E_n(H)=N_n=\lceil\al n\rceil$. If we let $E^H_n=E_n(H)-E_n(H-)$ then
\[
\hat E^H_n:=\frac{E^H_n-n\la^H}{\sqrt{n}}=\frac{N_n-E_n(H-)-n\al+np^{-1}\mu H}{\sqrt{n}}
=-\hat E_n(H-)+\frac{\lceil n\al\rceil-n\al}{\sqrt n}.
\]
Next, recall $\Xi(k)=\sum_{i=1}^k\xi_i$, $k\in\Z_+$. Let
\begin{equation}\label{58}
\hat\Xi_n(t)=\frac{\sum_{i=1}^{[nt]}(\xi_i-p)}{\sqrt{n}}
=\frac{\Xi([nt])-p[nt]}{\sqrt{n}},\qquad t\in\R_+.
\end{equation}
Now, $A_n(t)=\Xi\circ E_n(t)$ by \eqref{eq:arrival}. Let us consider this process
on $[0,H)$ separately from its jump at $H$.
For $t\in[0,H)$, use the above to write
\begin{align}\label{59}
\notag
\hat A_n(t)&:=\frac{A_n(t)-n\mu t}{\sqrt{n}}\\
\notag
&=\frac{\Xi(E_n(t))-pE_n(t)}{\sqrt{n}}+p\frac{E_n(t)-np^{-1}\mu t}{\sqrt{n}}\\
\notag
&=\frac{\Xi([n\bar E_n(t)])-p[n\bar E_n(t)]}{\sqrt{n}}+p\frac{E_n(t)-np^{-1}\mu t}{\sqrt{n}}\\
&=\hat\Xi_n(\bar E_n(t))+p\hat E_n(t),\hspace{12em} t\in[0,H).
\end{align}

The number of show ups at the time $H$ is
$A^H_n:=A_n(H)-A_n(H-)=\Xi(N_n)-\Xi(E_n(H-))$. Hence
\begin{align}\label{60}
\notag
\hat A^H_n&:=\frac{A^H_n-np\la^H}{\sqrt{n}}\\
\notag
&=\frac{\Xi(N_n)-np\al}{\sqrt{n}}-\frac{\Xi(E_n(H-))-pE_n(H-)}{\sqrt{n}}
-p\frac{E_n(H-)-np^{-1}\mu H}{\sqrt{n}}\\
&=\hat\Xi_n(\al)-\hat\Xi_n(\bar E_n(H-))-p\hat E_n(H-).
\end{align}
Next we let
\[
\hat S_n(t)=\frac{S_n(t)-n\mu t}{\sqrt n},\qquad t\ge0.
\]
We define the diffusion scale queue length, for $t\in[0,H)$ only, as
\[
\hat Q_n(t)=\frac{Q_n(t)}{\sqrt n},\qquad t\in[0,H).
\]
It is possible to consider the diffusion scale queue length for $t\ge H$
by first centering about $q^*(t)$ (as defined in Lemma~\ref{lem:upper-bound-convergence}) and then rescaling, but to avoid confusion,
we do not extend the process $\hat Q^n$ beyond the interval $[0,H)$.
We denote $q_n=n^{-1/2}(Q_n(H)-q^*(H))$.

By \eqref{eq:queue-length}, letting $I_n(t)=t-B_n(t)$ denote the cumulative idleness
process,
\begin{align}\label{61}
\notag
  \hat Q_n(t)&=n^{-1/2}(A_n(t)-S_n\circ B_n(t))\\
\notag
  &=\hat A_n(t)-\hat S_n(B_n(t))+n^{1/2}\mu I_n(t)\\
  &=\hat\Xi_n(\bar E_n(t))+p\hat E_n(t)-\hat S_n(B_n(t))+n^{1/2}\mu I_n(t),
  \qquad t\in[0,H),\\ \label{61+}
  &=\Gamma_1[p\hat E_n+X_n](t),\qquad t\in[0,H),
\end{align}
where \eqref{59} is used, and one denotes
\[
X_n(t)=\hat\Xi_n(\bar E_n(t))-\hat S_n(B_n(t)).
\]
The queue length dynamics for the BOP will later be derived from the above relation.

As for $q_n$, we can write it, using \eqref{60} and \eqref{61}, as
\begin{align}
\notag
q_n&=\hat Q_n(H-)+\hat A_n^H\\
\label{94}
&=-\hat S_n(B_n(H-))+\hat\Xi_n(\al)+n^{1/2}\mu I_n(H-)\\
\label{95}
&=\Gamma_1[p\hat E_n+X_n](H-)+\hat\Xi_n(\al)-\hat\Xi_n(\bar E_n(H-))-p\hat E_n(H-).
\end{align}

Next we develop equations for the two ingredients of the cost,
namely the overtime $[\tau_n-H]^+$ and the makespan, suitably
normalized at the diffusion scale.

For $t\ge H$, \eqref{eq:queue-length} is still valid, but $A_n(t)$ is simply
given by $A_n(H)$ since there are no arrivals after time $H$.
Moreover, the server is busy continuously on $[H,\tau_n)$, on the event $\tau_n>H$.
Hence $B_n(t)=B_n(H)+t-H$ for $t\in[H,\tau_n)$.
Clearly, $Q_n(H)=Q_n(H-)+A^H_n$.
For $t\in[H,\tau_n]$, the queue length is given by
\begin{align*}
Q_n(t)&=Q_n(H)-D_n(t)+D_n(H)\\
&=Q_n(H)-S_n(B_n(t))+S_n(B_n(H))\\
&=np\la^H+[Q_n(H-)+A^H_n-np\la^H]\\
&\quad -n\mu(t-H)-[S_n(B_n(t))-n\mu(B_n(t))]\\
&\quad +[S_n(B_n(H))-n\mu B_n(H)].
\end{align*}
Dividing by $\sqrt{n}$, for $t\in[H,\tau_n]$,
\begin{align}\label{63}
n^{-1/2}Q_n(t)&=n^{1/2}p\la^H+q_n-n^{1/2}\mu(t-H)
-\hat S_n(B_n(t))+\hat S_n(B_n(H)).
\end{align}
We have for $\tau_n$ the equation $Q_n(\tau_n)=0$
and for $\bar\tau$, $q^*(\bar\tau)=0$,
where we recall that for $t>H$, $q^*(t)=(p\la^H-\mu (t-H))\vee0
= (p\al-\mu H-\mu (t-H))\vee0=(p\al-\mu t)\vee 0$.
Hence $p\la^H-\mu(\bar\tau-H)=0$. Using these two relations in \eqref{63} gives
\[
0=n^{1/2}[-\mu(\tau_n-\bar\tau)]+q_n
-\hat S_n(B_n(\tau_n))+\hat S_n(B_n(H)).
\]
If we set
\begin{equation}\label{78}
\hat\tau_n=n^{1/2}(\tau_n-\bar\tau),
\end{equation}
then
\begin{equation}\label{62}
\mu\hat\tau_n=q_n-\hat S_n(B_n(\tau_n))+\hat S_n(B_n(H)).
\end{equation}
Equation \eqref{62} will be used to propose a formal limit of $\hat\tau_n$,
and later, to analyze rigorously the weak limit thereof.

Next, the FOP quantity for the makespan is
$\bar W=\int_H^{\bar\tau}(p\la^H-\mu(t-H))dt$,
where $\bar\tau=H+p\mu^{-1}\la^H$. Thus
\begin{align}
\notag
  \hat W_n&:=n^{1/2}(\bar W_n-\bar W)\\
\notag
  &=n^{-1/2}W_n-n^{1/2}\bar W\\
\label{65}
  &=\hat W_n(1)+\hat W_n(2),\\
\label{66}
  \hat W_n(1)&:=\int_0^{H\wedge\tau_n}\hat Q_n(t)dt,\\
\label{92}
  \hat W_n(2)&:=
  n^{-1/2}\int_{H\wedge\tau_n}^{\tau_n}Q_n(t)dt-n^{1/2}\int_H^{\bar\tau}(p\la^H-\mu(t-H))dt.
\end{align}
On the event $\tau_n<H$, the first term in the expression for $\hat W_n(2)$
is zero, and we obtain $\hat W_n(2)=-n^{1/2}\bar W$.
Next, consider the event $\tau_n\ge H$. Then
\begin{align}\label{90}
  \hat W_n(2)&=\int_H^{\tau_n}[n^{-1/2}Q_n(t)-n^{1/2}(p\la^H-\mu(t-H))]dt
  +n^{1/2}\int_{\bar\tau}^{\tau_n}(p\la^H-\mu(t-H))dt.
\end{align}
A use of \eqref{63} and the computed value of $\bar\tau$ gives
\begin{equation}\label{64}
\hat W_n(2)=\int_H^{\tau_n}[q_n-\hat S_n(B_n(t))
+\hat S_n(B_n(H))]dt
-\frac{\mu}{2}n^{1/2}(\tau_n-\bar\tau)^2.
\end{equation}

To derive the BOP, note first that the FCLT applies to the processes
$\hat\Xi_n$ and $\hat S_n$. That is,
let $X^{(1)}$ and $X^{(2)}$ be mutually independent 1-dimensional BMs
with zero drift and diffusion coefficient $(p(1-p))^{1/2}$ and $\mu^{1/2}C_S$, respectively.
Then by the FCLT, $(\hat\Xi_n,\hat S_n)\To(X^{(1)},X^{(2)})$ \cite[\S 17]{Bi68}.

We take {\it formal} limits in equation \eqref{61}. Denote by $Q$, $U$ and $L$
limits of the processes $\hat Q_n$, $p\hat E_n$ and $n^{1/2}\mu I_n$, respectively,
on the time interval $[0,H)$. Denote $\tilde\mu=p^{-1}\mu$,
and {\em approximate} $\bar E_n(t)$ as $\tilde\mu t$,
and $B_n(t)$ as $t$, $t\in[0,H)$. Then, we expect the following relationship to hold in the limit
\begin{equation}\label{69}
Q(t)=U(t)+X^{(1)}(\tilde\mu t)-X^{(2)}(t)+L(t), \qquad t\in[0,H).
\end{equation}
Moreover,
\begin{equation}
  \label{67}
  Q(t)\ge0,\ t\in[0,H), \text{ and } \int_{[0,H)}Q(t)dL(t)=0.
\end{equation}
From \eqref{60}, letting $A^H$ be a weak limit of $\hat A^H_n$,
\[
A^H=X^{(1)}(\al)-X^{(1)}(\tilde\mu H)-U(H-).
\]
To obtain an expression for $\hat\tau$, a weak limit of $\hat\tau_n$,
use \eqref{62} to write
\begin{align*}
\hat\tau&=\mu^{-1}(Q(H-)+A^H-X^{(2)}(\bar\tau)+X^{(2)}(H))\\
&=\mu^{-1}(X^{(1)}(\al)-X^{(2)}(\bar\tau)+L(H-)).
\end{align*}
With $\hat W$, $\hat W(1)$ and $\hat W(2)$ representing limits of $\hat W_n$,
$\hat W_n(1)$ and $\hat W_n(2)$, respectively, we have from \eqref{65},
$\hat W(1)=\int_0^HQ(t)dt$ (since $\tau_n$ converges to $\bar\tau\ge H$),
while by \eqref{64},
\[
\hat W(2)=\int_H^{\bar\tau}(Q(H-)+\hat A^H-X^{(2)}(t)+X^{(2)}(H))dt
=\int_H^{\bar\tau}(X^{(1)}(1)-X^{(2)}(t)+L(H-))dt.
\]
The cost is thus given by
\begin{align}\label{68}
\notag
  \hat J(U)&=c_w\E[\hat W(1)+\hat W(2)]+c_o\E[\hat\tau]\\
\notag
  &=c_w\E\int_0^HQ(t)dt+c_w(\bar\tau-H)\E[L(H-)]+c_o\mu^{-1}\E[L(H-)]\\
  &=c_w\E\int_0^HQ(t)dt+\tilde c_o\E[L(H-)],
\end{align}
where $\tilde c_o=c_w(\bar\tau-H)+c_o\mu^{-1}$.
Finally, we can simplify \eqref{69} by considering $X$, a BM with drift zero
and diffusion coefficient
$\sigma=(\tilde\mu p(1-p)+\mu C_S^2)^{1/2}=\mu^{1/2}(1-p+C_S^2)^{1/2}$
in place of the two BM terms, and write
\begin{equation}
  \label{70}
  Q(t)=U(t)+X(t)+L(t),\qquad t\in[0,H).
\end{equation}
We thus let $\calU$ denote the collection of RCLL functions $u:[0,H)\to\R$,
and note that given $u\in\calU$ \eqref{70} and \eqref{67} uniquely define $Q$ and $L$ in terms
of $X$.

We can now state the BOP of interest, as a problem involving
\eqref{67}, \eqref{68} and \eqref{70}, with value given by
\begin{equation}
  \label{71}
  \hat V=\inf_{U\in\calU}\hat J(U).
\end{equation}

{\it Remark.}
Instead of $\calU$ being functions defined on $[0,H)$ we can work with
functions on $[0,H]$ and replace the term $\E[L(H-)]$ by $\E[L(H)]$
in \eqref{68}. This will not change the value $\hat V$,
due to the fact that having a jump $L(H)-L(H-)>0$ can only
increase the cost $\hat J$ as compared to having $L(H)=L(H-)$
(the jump cannot be negative since $L$ is non-decreasing).
Throughout what follows we shall work with this slightly modified
definition of $\calU$ and $\hat J(U)$.

{\it Remark.}
We can present the optimization problem in a way that the cost is more explicit,
and moreover makes it easy to see that it is a convex optimization problem.
The pair of equations \eqref{67} and \eqref{70} is related to the Skorohod Problem
on the half line. Namely, $Q=\Gamma_1[U+X]$. Thus
\begin{equation}
  \label{72}
  Q(t)=U(t)+X(t)-\inf_{s\in[0,t]}[(U(s)+X(s))\wedge 0],
  \qquad
  L(t)=Q(t)-U(t)-X(t).
\end{equation}
We can therefore write $\hat J$ as
\[
\hat J(U)=c_w\E\int_0^H\Gamma_1[U+X](t)dt+\tilde c_o\E[\Gamma_1[U+X](H)-U(H)-X(H)].
\]

\subsection{Large time solution of the BOP}\label{sec42}

In this subsection we analyze the BOP at the large $H$ limit.
Note carefully that Assumption \eqref{eq:overload} puts a restriction on $H$,
namely $p^{-1}\mu H < \alpha$.
Thus $H$ cannot be taken arbitrarily large without modifying
$(\mu,p,\alpha)$. In our treatment $\mu$ and $p$ remain fixed
and $\alpha$ and $H$ grow so that the assumption remains valid.
However, this issue is not significant in this subsection where we only
work with the BOP itself, because the parameter $\alpha$
does not show up in it. It does become relevant in later sections.

Recall that $X$ is a $(0,\sig)$-BM,
$\calU_H$ denotes the collection of RCLL functions
$[0,H]\to\R$, and for $U\in\calU_H$ let $L=L^U$ and $Q=Q^U$ be defined as
\[
L_t=\sup_{s\in[0,t]}(-X(s)-U(s))^+,\qquad Q(t)=X(t)+U(t)+L(t),\qquad t\in[0,H].
\]
Let also
\begin{equation}\label{75}
\hat J_H(U)= \frac{c_w}{H}\E\int_0^HQ(t)dt+\frac{\tilde c_o}{H}\E[Q(H)-U(H)],\qquad U\in\calU_H,
\end{equation}
\[
\hat V_H=\inf_{U\in\calU_H}\hat J_H(U).
\]
Denote by $\calU_H^{\rm lin}$ the collection of linear functions $U(t)=\beta t$, $t\in[0,H]$
for some $\beta\in\R$. Note that the process $Q$ corresponding to such a control
is a reflected BM with drift $\beta$ and diffusion coefficient $\sigma$.

For $\beta<0$ let $m_{{\rm RBM}(\beta)}(dx)=-\frac{2\beta}{\sigma^2}e^{2\beta x/\sig^2}dx$,
$x\in[0,\iy)$. This probability measure on $[0,\iy)$
is the stationary distribution of RBM with
drift $\beta<0$ and diffusion coefficient $\sigma$. Let
\begin{equation}\label{30}
V^*=\inf_{\beta<0} \left[c_w\int x\,m_{{\rm RBM}(\beta)}(dx)-\tilde
  c_o\beta \right].
\end{equation}
We next establish that for the large horizon BOP it is sufficient to consider control functions in $\calU_H^{\rm lin}$.
\begin{proposition}
  \label{prop4}
One has
\begin{align}
\lim_{H\to\iy}\hat V_H
\label{40}
&=\lim_{H\to\iy}\inf_{U\in\calU_H^{\rm lin}}\hat J_H(U)
\\
\label{41}
&=V^*=\sig\sqrt{2c_w\tilde c_o}.
\end{align}
Moreover,
$\beta^*=-\sig\sqrt{c_w/(2\tilde c_o)}$ is optimal for both the expressions in \eqref{30} and \eqref{40},
that is, with $U^*(t)=\beta^*t$, $\lim_{H\to\iy}\hat J_H(U^*)=V^*$
and $c_w\int x\,m_{{\rm RBM}(\beta^*)}(dx)-\tilde c_o\beta^*=V^*$.
\end{proposition}

The proof is based on several lemmas.
The first is concerned with large time behavior of RBM and computes $V^*$.
\begin{lemma}
  \label{lem4--}
1. For each $\beta\in\R$, let $Q^\beta_t$ be a $(\beta,\sig)$-RBM starting at the origin.
Then
\[
\lim_{t\to\iy}\inf_{\beta\in\R}[c_w\E[Q^\beta_t]-\tilde c_o\beta]=V^*
\]
2. The infimum in \eqref{30} is attained at $\beta^*=-\sig\sqrt{c_w/(2\tilde c_o)}$ and
$V^*=\sig\sqrt{2c_w\tilde c_o}$.
\end{lemma}

The following lemma shows that one can focus on controls under which
$\E[Q_H]$ is sublinear in $H$. Stated precisely, we show
\begin{lemma}\label{lem4-}
For every $\eps>0$ there exists $H_0$ such that for $H>H_0$
and $U$ for which $\E[Q^U_H]\ge\eps H$
one can find $\tilde U$ with
$\E[Q^{\tilde U}_H]<\eps H$ and $\hat J_H(\tilde U)\le\hat J_H(U)+e_1(H)$.
Here, $e_1(H)$ does not depend on $U$ and converges to zero as $H\to\infty$.
\end{lemma}

The following lemma argues that one may focus on controls that are constant
on an initial and terminal intervals. More precisely,
given $z,k,H\in(0,\iy)$, $2z<H$, define the class of controls
$\calU^\#(z,k,H)$ as the collection of members $U\in\calU_H$ which satisfy
$U_t=k$ for $t\in[0,z)$ and $U_t=U_H$ for $t\in[H-z,H]$.
\begin{lemma}\label{lem40}
Given $z,k,H$ and $U\in\calU_H$, let $U^\#\in\calU^\#(z,k,H)$ be defined as
\[
U^\#_t=
\begin{cases}
  k,& t\in[0,z),\\ U_t,& t\in[z,H-z),\\ U_H,& t\in[H-z,H].
\end{cases}
\]
Fix $\eps>0$, let $H_0=H_0(\eps)$ be as in Lemma \ref{lem4-}, and consider $H>H_0$
and $U\in\calU_H$ for which $\E[Q_H]\le\eps H$. Then
\begin{equation}
  \label{76}
  \hat J_H(U^\#)\le\hat J_H(U)+e_2(z,k,H),
\end{equation}
where
\begin{equation}\label{77}
\limsup_He_2(z,k,H)\le c_w\E[r]+c_w\eps z,
\end{equation}
\[
r=r(k,z,X)=\sup_{s\in[0,z)}(-k-X_s)^+.
\]
\end{lemma}

The following lemma relates the large time behavior of $\hat J_H(U)$
for $U$ as in Lemma \ref{lem40} to the expression $V^*$.
\begin{lemma}
  \label{lem4}
One has
\[
\liminf_{H\to\iy}\inf_{U\in\calU^\#(z,k,H)}\hat J_H(U)\ge V^*-e_3(z),
\]
where $e_3(z)\to0$ as $z\to\iy$.
\end{lemma}

\skp

\noi{\bf Proof of Proposition \ref{prop4}.}
Using Lemmas \ref{lem4-}, \ref{lem40} and \ref{lem4}, for any $\eps>0$, $z>0$ and $k>0$,
\begin{align*}
  \liminf_{H\to\iy}\hat V_H&\ge\liminf_{H\to\iy}\inf_{U\in\calU_H:EQ^U_H<\eps H}\hat J_H(U)\\
  &\ge\liminf_{H\to\iy}\inf_{U\in\calU^\#(z,k,H):EQ^U_H<\eps H}\hat J_H(U)-\limsup_{H\to\iy}e_2(z,k,H)\\
  &\ge V^*-\limsup_{H\to\iy}e_2(z,k,H)-e_3(z),
\end{align*}
for $e_2$ and $e_3$ as in these lemmas.
Thus
\begin{align*}
  \liminf_{H\to\iy}\hat V_H
  &\ge V^*-c_wE\Big[\sup_{s\in[0,z)}(-k-X_s)^+\Big]-c_w\eps z-e_3(z).
\end{align*}
We refer to the last three terms on the RHS above as the first, second and third error terms
in the order at which they appear.
We first take $\eps\to0$ so that the second error term vanishes. Then we take $k\to\iy$
to have the first error term vanish, as a direct consequence of
$\sup_{s\in[0,z)}(-k-X_s)^+\le(-k+\|X\|_z)^+$ and $E[\|X\|_z]<\iy$.
Finally we take $z\to\iy$ and the third term vanishes.
We have thus shown that $\liminf_{H\to\iy}\hat V_H\ge V^*$. For a matching upper bound
we simply select $U\in\calU^{\rm lin}_H$ with $\beta=\beta^*$ and appeal to Lemma \ref{lem4--}.
We conclude that $\lim_{H\to\iy}\hat V_H=V^*$.
\hfill$\Box$

\skp

\noi{\bf Proof of Lemma \ref{lem4--}.}
For each $t\ge0$, The CDF of $Q^\beta_t$ is given by
\[
P(Q^\beta_t\le y)=\Ph\Big(\frac{y-\beta t}{\sig t^{1/2}}\Big)
-e^{2\beta y/\sig^2}\Ph\Big(\frac{-y-\beta t}{\sig t^{1/2}}\Big),
\qquad y\ge0,
\]
where $\Ph$ is the standard normal CDF (\cite{Ha1985} page 15).
For fixed $\beta<0$,
the limit distribution as $t\to\iy$ is exponential with mean $\sig^2/(2|\beta|)$.
Moreover, it can be directly checked that the CDF is monotone decreasing
in $t$. Hence by monotone convergence,
the expectation $\E[Q^\beta_t]$ converges, as $t\to\iy$, to $\sig^2/(2|\beta|)$,
provided $\beta<0$. For $\beta\ge0$, $P(Q^\beta_t\le y)\to0$ for all $y$, hence
$\E[Q^\beta_t]\to\iy$. Thus, denoting $F(t,\beta)=c_w\E[Q^\beta_t]-\tilde c_o\beta$
and $F(\beta)=c_w\sig^2/(2|\beta|)+\tilde c_o|\beta|$ for $\beta<0$, $F(\beta)=\iy$
for $\beta\ge0$,
we have the pointwise convergence $\lim_{t\to\iy}F(t,\beta)=F(\beta)$.

Our goal now is to show
\begin{equation}\label{31}
\lim_{t\to\iy}\inf_{\beta\in\R}F(t,\beta)=\inf_{\beta<0}F(\beta).
\end{equation}
We achieve this in three steps. First we show that $\inf_{\beta\ge0}F(t,\beta)\to\iy$ as $t\to\iy$.
Then we argue that there exist $-\iy<a<-1<b<0$
such that, for all large $t$, $\inf_{\beta\in\R}F(t,\beta)=\inf_{\beta\in[a,b]}F(t,\beta)$.
Then we are in a position to use Dini's theorem
to argue that the order of the $t$-limit and the $\beta$-infimum
can be interchanged. We shall use an additional monotonicity property. It can be readily checked by
the above CDF formula that $\beta\to P(Q^\beta_t\le y)$ is monotone decreasing (for each
$t$ and $y$). Hence $\beta\to E[Q^\beta_t]$ is monotone increasing (for each $t$).

For the first step alluded to above,
the pointwise convergence $EQ^0_t\to\iy$ as $t\to\iy$ can be used to deduce
$\inf_{\beta\in[0,1]}[c_wEQ^\beta_t-\tilde c_o\beta]\to\iy$ because for $\beta\in[0,1]$
and all $t$ we have $c_wEQ^\beta_t-\tilde c_o\beta\ge c_wEQ^0_t-\tilde c_o\to\iy$ as $t\to\iy$.
To show $\inf_{\beta\in(1,\iy)}[c_wEQ^\beta_t-\tilde c_o\beta]\to\iy$ we
argue as follows. By the formula for $Q^\beta_t$ we have the lower bound $Q^\beta_t\ge X_t+\beta t$.
Hence for all $\beta>1$, $c_wQ^\beta_t-\tilde c_o\beta\ge c_wX_t+\beta(c_wt-\tilde c_o)$. Hence for $t>\tilde c_o/c_w$,
$c_wQ^\beta_t-\tilde c_o\beta\ge c_wX_t+(c_wt-\tilde c_o)$. Taking expectation gives
$\inf_{\beta\in(-\iy,-1)}[c_wEQ^\beta_t-\tilde c_o\beta]\ge c_wt-\tilde c_o\to\iy$ as $t\to\iy$.

For the next step, fix $u>0$ and let $-\iy<a<-1<b<0$ such that
$\tilde c_o|a|>u$, $\tilde c_o|b|<1$ and $F(b)>u$ (note that $F(0-)=\iy$).
Then for any $\beta<a$, $F(t,\beta)\ge \tilde c_o|\beta|>u$. Next, consider $\beta\in(b,0)$.
The pointwise convergence of $F(t,b)$ to $F(b)$ implies that for some $t_0$
and all $t\ge t_0$, $F(t,b)>u-1$ and, since $\tilde c_o|b|<1$, $c_wEQ^b_t>u-2$.
Using again the monotonicity in $\beta$,
$F(t,\beta)\ge c_wEQ^\beta_t>u-2$ for $\beta\in(b,0)$.

Since for fixed $\beta<0$ the limit $\lim_{t\to\iy}F(t,\beta)$ is finite,
and the constant $u$ is arbitrary,
it follows that $a$ and $b$ as above can be found so that the infimum is achieved
in $[a,b]$ for all large $t$.

Next, $F(t,\beta)$ is continuous in $\beta$ for $\beta$
in the compact interval $[a,b]$, and
these functions converge pointwise as $t\to\iy$ to a continuous function
$F(\beta)$. Hence by Dini's theorem, the convergence is uniform in $\beta$.
We conclude that,
as $t\to\iy$, $\inf_{\beta\in[a,b]}F(t,\beta)\to\inf_{\beta\in[a,b]}F(\beta)$.
This proves part 1 of the lemma.

It remains to solve the optimization problem \eqref{30},
or equivalently the RHS of \eqref{31}.
As already stated, for $\beta<0$, the first moment of $m_{{\rm RBM}(\beta)}$ is given by
$\sig^2/(2|\beta|)$. We are therefore interested in minimizing
\[
c_w\frac{\sig^2}{2|\beta|}+\tilde c_o|\beta|
\]
over $\beta\in(-\iy,0)$.
By a direct calculation, the minimum is attained at $\beta^*=-\sig(c_w/2\tilde c_o)^{1/2}$
and is given by $\sig(2c_w\tilde c_o)^{1/2}$.
\qed

\skp

\noi{\bf Proof of Lemma \ref{lem4-}.}
Fix $\eps>0$. Given any $H$ and any $U\in\calU_H$
we have for $Q=Q^U$ the relation
\[
Q_t=U_t+X_t+\sup_{s\le t}(-U_s-X_s)^+,\qquad t\in[0,H].
\]
Define $q=q^U$ as
\[
q_t=U_t+\sup_{s\le t}(-U_s)^+,\qquad t\in[0,H].
\]
Then
\begin{equation}\label{74}
\|Q-q\|_H\le2\|X\|_H.
\end{equation}
Consider $U$ for which $\E[Q_H]\ge\eps H$.
Then $q_H\ge\eps H-2\E[\|X\|_H]$.
Fix $H_0$ so large that for every $H>H_0$, $2\E[\|X\|_H]<\eps H/4$,
and consider in what follows only $H>H_0$.
Then $q_H>3\eps H/4$.
Let $\tilde U$ be defined as $\tilde U=U$ on $[0,H)$ and $\tilde U_H=U_H-q_H$
(here, $q=q^U$).
Denote $\tilde Q=Q^{\tilde U}$ and $\tilde q=q^{\tilde U}$.
Clearly, on $[0,H)$, we have $\tilde q=q$ and $\tilde Q=Q$.
As for the time $H$, we have
\begin{align*}
  \tilde q_H&=U_H-q_H+\sup_{s\le H}(-U_s+q_H1_{\{s=H\}})^+.
\end{align*}
As a result,
\[
\tilde Q_H=\tilde q_H+\tilde Q_H-\tilde q_H=\tilde Q_H-\tilde q_H\le 2\|X\|_H,
\]
by \eqref{74}. This shows $\E[\tilde Q_H]\le2\E[\|X\|_H]<\eps H/4$.
Moreover, by \eqref{75},
\[
\hat J_H(\tilde U)\le\hat J_H(U)+\frac{\tilde c_o}{H}2\E[\|X\|_H]\le
\hat J_H(U)+c_4H^{-1/2},
\]
for a suitable constant $c_4$. This proves the lemma.
\hfill$\Box$

\skp

\noi{\bf Proof of Lemma \ref{lem40}.}
Denote $Q=Q^U$ and $Q^\#=Q^{U^\#}$.
First, we provide lower estimates on $Q-Q^\#$
on each of the three intervals separately.

{\it The interval $[0,z)$.} Here we use the trivial lower bound $Q_t\ge0$. As for $Q^\#$,
\[
Q^\#_t=k+X_t+\sup_{s\le t}(-k-X_s)^+\le2k+2\|X\|_z.
\]

{\it The interval $[z,H-z)$.}
We have
\begin{align*}
  Q_t&=U_t+X_t+\sup_{s\le t}(-U_s-X_s)^+\\
  &\ge U_t+X_t+\sup_{s\in[z,t]}(-U_s-X_s)^+\\
  &\ge U_t+X_t+\max\Big[r,\sup_{s\in[z,t]}(-U_s-X_s)^+\Big]-r,
\end{align*}
where we have used
the fact that $a\ge a\vee b-b$ provided $a\ge0$ and $b\ge0$.
Recalling that for $t\in[z,H-z)$ $U$ and $U^\#$ agree and that $U^\#=k$
on $[0,z)$, the expression above is equal to
\[
U^\#_t+X_t+\sup_{s\in[0,t]}(-U^\#_s-X_s)^+-r.
\]
It follows that
\[
Q_t\ge Q^\#_t-r.
\]

{\it The interval $[H-z,H]$.}
In fact, we will only need a lower estimate on $Q_H-Q^\#_H$.
We have
\[
Q_H=U_H+X_H+\sup_{s\le H}(-U_s-X_s)^+,
\]
while
\[
Q^\#_H=U_H+X_H+\max\Big[\sup_{s<H-z}(-U^\#_s-X_s)^+,\sup_{s\in[H-z,H]}(-U_H-X_s)^+\Big].
\]
Hence
\begin{align*}
  Q_H&\ge U_H+X_H+\sup_{s\in[z,H]}(-U_s-X_s)^+\\
  &\ge U_H+X_H+\max\Big[r,\sup_{s\in[z,H]}(-U_s-X_s)^+\Big]-r\\
  &\ge U_H+X_H
  +\max\Big[r,\sup_{s\in[z,H-z)}(-U_s-X_s)^+,\sup_{s\in[H-z,H]}(-U_H-X_s)^+\Big]-r-\hat r,
\end{align*}
where $\hat r=\sup_{s\in[H-z,H]}|X_s-X_H|$.
This shows
\[
Q_H\ge Q^\#_H-r-\hat r.
\]

Next,
\begin{align*}
  \hat J_H(U)&= \frac{c_w}{H}\E\int_0^HQ_sds+\frac{\tilde c_o}{H}\E[Q_H-U_H]\\
  &\ge\frac{c_w}{H}\E\int_z^{H-z}Q_sds+\frac{\tilde c_o}{H}\E[Q_H-U_H]\\
  &\ge \frac{c_w}{H}\Big[\E\int_z^{H-z}Q^\#_sds-(H-2z)\E[r]\Big]
  +\frac{\tilde c_o}{H}\{\E[Q^\#_H-U^\#_H]-\E[r+\hat r]\}.
\end{align*}
Also,
\begin{align*}
  \hat J_H(U^\#)&= \frac{c_w}{H}\E\int_0^HQ^\#_sds+\frac{\tilde c_o}{H}\E[Q^\#_H-U^\#_H]\\
  &\le \frac{c_w}{H}\Big[2k+2\E[\|X\|_z]+\E\int_z^{H-z}Q^\#_sds+z\E[Q^\#_H+\hat r]\Big]
  +\frac{\tilde c_o}{H}\E[Q^\#_H-U^\#_H],
\end{align*}
where we used that for $t\in[H-z,H]$ one has $Q^\#_t=X_t+U^\#_t+L^\#_t
\le X_t+U^\#_H+L^\#_H=X_t-X_H+Q^\#_H$, whence $Q^\#_t\le Q^\#_H+\hat r$.
Combine these two bounds to obtain
\begin{align*}
\hat J_H(U^\#)-\hat J_H(U)
&\le\frac{c_w}{H}\Big[(H-2z)\E[r]+2k+2\E[\|X\|_z]+z\E[Q^\#_H+\hat r]\Big]+\frac{\tilde c_o}{H}\E[r+\hat r]\\
&\le\frac{c_w}{H}\Big[(H-2z)\E[r]+2k+2\E[\|X\|_z]+z\E[Q_H+r+2\hat r]\Big]+\frac{\tilde c_o}{H}\E[r+\hat r]\\
&\le\frac{c_w}{H}\Big[(H-2z)\E[r]+2k+2\E[\|X\|_z]+z\eps H+z\E[r+2\hat r]\Big]+\frac{\tilde c_o}{H}\E[r+\hat r].
\end{align*}
Denote by $e_2(z,k,H)$ the expression on the last line above. Note that while $\hat r$
depends on $H$, its expectation does not (and is finite). Thus $e_2$ satisfies \eqref{77}.
\hfill$\Box$

\skp

\noindent{\bf Proof of Lemma \ref{lem4}.}
Fix $z,k,H$ and a control $U\in\calU^\#(z,k,H)$.
Then $Q_t=Q^U_t=X_t+U_t+\sup_{s\le t}[-X_s-U_s]^+\ge\sup_{s\le t}[X_t-X_s+U_t-U_s]$.
Hence, for $z+t\le H$,
\begin{align*}
\int_0^{z+t}Q_sds
&\ge \int_z^{z+t}Q_sds
\\
&\ge\int_z^{z+t}\sup_{\theta\in[0,s]}[X_s-X_{s-\theta}+U_s-U_{s-\theta}]ds
\\
&\ge\int_z^{z+t}\sup_{\theta\in[0,z]}[X_s-X_{s-\theta}+U_s-U_{s-\theta}]ds.
\end{align*}
For each $s\ge z$, the stochastic process $\{X_s-X_{s-\theta}\}_{\theta\in[0,z]}$
is equal in law to $\{X_z-X_{z-\theta}\}_{\theta\in[0,z]}$.
As a result,
\[
\E\int_0^{z+t}Q_sds
\ge \E\int_z^{z+t}\sup_{\theta\in[0,z]}[X_z-X_{z-\theta}+U_s-U_{s-\theta}]ds.
\]
We now use the inequality
\[
\int_a^b \sup_\theta f(s,\theta)ds
\ge\sup_\theta\int_a^b f(s,\theta)ds.
\]
This gives
\begin{align*}
\frac{1}{t}\E\int_0^{z+t}Q_sds
&\ge \frac{1}{t}\E\sup_{\theta\in[0,z]}
\int_z^{z+t}[X_z-X_{z-\theta}+U_s-U_{s-\theta}]ds
\\
&=\E\sup_{\theta\in[0,z]}
\Big\{[X_z-X_{z-\theta}]+\frac{1}{t}\int_z^{z+t}[U_s-U_{s-\theta}]ds\Big\}.
\end{align*}
We have
\[
\int_z^{z+t}[U_s-U_{s-\theta}]ds=\int_z^{z+t}U_sds
- \int_{z-\theta}^{z-\theta+t}U_sds
=\int_{z-\theta+t}^{z+t}U_sds-\int_{z-\theta}^z U_sds.
\]
Moreover,
by the assumption on $U$, we have $U=k$ on $[0,z)$, $U=U_H$ on $[H-z,H]$.
Thus with $t+z=H$,
\[
\frac{1}{H-z}\E\int_0^HQ_sds\ge
\E\sup_{\theta\in[0,z]}
\Big\{[X_z-X_{z-\theta}]+\theta\frac{1}{H-z}(U_H-k)\Big\}.
\]
Thus
\begin{align*}
\hat J_H(U)&= \frac{c_w}H \E\int_0^HQ_tdt+\frac{\tilde c_o}{H}\E[Q_H-U_H]
\\
&\ge c_w\E\sup_{\theta\in[0,z]}
\Big\{\frac{H-z}{H}[X_z-X_{z-\theta}]+\theta(\frac{U_H}{H}-\frac{k}{H})\Big\}
-\tilde c_o\frac{U_H}{H}\\
&\ge\inf_{\beta\in\R}
\Big\{c_w\E\sup_{\theta\in[0,z]}
\Big\{\frac{H-z}{H}[X_z-X_{z-\theta}]+\theta(\beta-\frac{k}{H})\Big\}
-\tilde c_o\beta\Big\}.
\end{align*}
Denoting $\del=z/H$,
\begin{align*}
\hat J_H(U)
&\ge\inf_{\beta\in\R}
\Big\{c_w\E\sup_{\theta\in[0,z]}
\Big\{(1-\del)[X_z-X_{z-\theta}]+\theta\beta\Big\}
-\tilde c_o\beta\Big\}-k\del\\
&\ge\inf_{\beta\in\R}
\Big\{c_w\E\sup_{\theta\in[0,z]}
\Big\{[X_z-X_{z-\theta}]+\theta\beta\Big\}
-\tilde c_o\beta\Big\}-k\del-2\del \E\|X\|_z.
\end{align*}
Send $H\to\iy$ (hence $\del\to0$) to obtain
\[
\liminf_{H\to\iy}\inf_{U\in\calU^\#(z,k,H)}\hat J_H(U)\ge
\inf_{\beta\in\R}\La(z,\beta),
\]
where
\[
\La(z,\beta)=
\Big\{c_w\E\sup_{\theta\in[0,z]}
\Big\{[X_z-X_{z-\theta}]+\theta\beta\Big\}
-\tilde c_o\beta\Big\}.
\]
Now, if we let $Q_z=\sup_{\theta\in[0,z]}\{[X_z-X_{z-\theta}]+\theta\beta\}$,
then $Q_z$ is a $(\beta,\sig)$-RBM starting from zero.
Hence by Lemma \ref{lem4--}(1), $\lim_{z\to\iy}\inf_{\beta\in\R}\La(z,\beta)=V^*$.
This proves the lemma.
\hfill$\Box$

\subsection{Lower bound on the diffusion scale cost}\label{sec43}

Recall from \eqref{75}
the definitions of $\hat J_H$ and $\hat V_H$, the cost and value
of the BOP. Also recall
the diffusion scale cost $\hat J_{n,H}(\{T_i\})$ and value
$\hat V_{n,H}=n^{1/2}[V_{n,H}-\bar V_H]$ where $V_{n,H}$ is defined
in \eqref{51} and $\bar V_H$ is the FOP value defined in~\eqref{eq:FOP}.
This subsection is devoted to proving the following result.
\begin{theorem}
  \label{prop5}
Let Assumption \ref{assn1} hold.
Fix $H$. Then
\begin{equation}\label{73}
H^{-1}\liminf_{n\to\iy}\hat V_{n,H}\ge\hat V_H.
\end{equation}
\end{theorem}

Assumption \ref{assn1} is in force throughout this subsection.
(It is used in the proof of Lemma \ref{lem6} below).

We say that a sequence $\{T_i^n\}\in\calT_n$, $n\in\N$
{\it achieves the limit inferior in \eqref{73}} if
\begin{equation}\label{80}
\lim_{n\to\iy}\hat J_{n,H}(\{T_i^n\})=\liminf_{n\to\iy}\hat V_{n,H}.
\end{equation}
If the expression on the LHS of \eqref{73}
is infinite then there is nothing to prove. Hence we may, and will,
assume that for any such $\{T_i^n\}$,
the sequence $\hat J_{n,H}(\{T_i^n\})$ of \eqref{80} is bounded.
The following lemma provides verious convergence results,
based to a large extent on the boundedness of the sequence of costs.

\begin{lemma}
  \label{lem5}
There exists a sequence $\{T^n_i\}\subset\calT_n$, $n\in\N$ that achieves
the limit inferior in \eqref{73} and for which assertions (i)--(iv) below hold.
\\
i. $\sup_n\E[(\hat\tau_n^-)^2]\vee\E[\hat\tau_n^+]<\iy$
(in particular, $\hat\tau_n$ are tight and $\tau_n\To\bar\tau$, as $n\to\iy$).
\\
ii. $\sup_{t\in[0,\tau_n]}|B_n(t)-t|\To0$, as $n\to\iy$.
\\
iii. $\sup_n\int_0^H\Gamma_1[p\hat E_n](t)dt<\iy$
and $\sup_n\{\Gamma_1[p\hat E_n](H-)-p\hat E_n(H-)\}<\iy$.
\\
iv. $\sup_{t\in[0,H)}|\bar E_n(t)-p^{-1}\mu t|\to0$ and $|\bar E_n^H-\la^H|\to0$,
as $n\to\iy$.
\end{lemma}

Fix a sequence $\{T_i^n\}$ as in Lemma \ref{lem5}.
Let all the processes and RVs such as $W_n$, $Q_n$, $\tau_n$,
etc.\ denote those associated with the schedule $\{T_i^n\}$, for each $n$.
Recall the diffusion scale expression
$\hat W_n$ from \eqref{65}. Also recall $\hat W_n=\hat W_n(1)+\hat W_n(2)$.
Then by \eqref{43},
\begin{equation}\label{79}
\hat J_{n,H}(\{T_i^n\})=\E[c_w\hat W_n(1)+c_w\hat W_n(2)+c_o\hat\tau_n].
\end{equation}
Recall that $(\hat\Xi_n,\hat S_n)\To(X^{(1)},X^{(2)})$, that these two BMs
are mutually independent, and that, by its definition, $X$ is equal in law to
the sum $X^{(1)}(\tilde\mu\cdot)+X^{(2)}(\cdot)$.

In the derivation of the BOP in \S \ref{sec41}, we took formal limits
in the equations that describe the scaled processes, such as \eqref{61}
for $\hat Q_n$, \eqref{62} for $\hat\tau_n$, etc.
{\it We are unable to turn this into a rigorous argument} in a straightforward manner,
because there is no apparent pre-compactness for the sequence of functions $\{\hat E_n\}$.
For example, in \eqref{61}, there is no justification to replace
the limit of the term $p\hat E_n$ by some control $U$.
We therefore take a different route, where we are able to provide
a lower bound in which the error term converges to zero as $n\to\iy$
due only to the convergence of the stochastic processes
and RVs involved (such as $\tau_n$, $\hat \Xi_n$), not relying on any convergence
associated with $\hat E_n$.

An outline of the argument is as follows. We
appeal to Skorohod's representation theorem and derive (in Lemma \ref{lem6} below)
a bound of the form
\begin{equation}\label{82}
\E[c_w\hat W_n(1)+c_w\hat W_n(2)+c_o\hat \tau_n]\ge H\hat J_H(p\hat E_n)-\E[\eps_n],
\end{equation}
where $\eps_n$ is a sequence of RVs satisfying
$\E[\eps_n]\to0$ as $n\to\iy$.
This estimate is based on the closeness of $(\hat\Xi_n,\hat S_n)$
to $(X^{(1)},X^{(2)})$ but not on $\hat E_n$ being close to any candidate limit.
Now, for each $n$, $p\hat E_n$ is a member of $\calU_H$. Thus, using \eqref{79},
it follows from the definition of $\hat V_H$ that
\begin{equation}\label{83}
\hat J_{n,H}(\{T_i^n\})\ge H\hat V_H-\E[\eps_n].
\end{equation}
In view of the fact that $\E[\eps_n]\to0$, the result will follow.

\skp

Toward stating Lemma \ref{lem6}, note that
the convergence $\hat\Xi_n\To X^{(1)}$ and the one stated in Lemma~\ref{lem5}(iv)
imply that $\hat\Xi_n\circ\bar E_n\To X^{(1)}(\tilde\mu\cdot)$ (recall that
$\tilde\mu=p^{-1}\mu$). By Lemma \ref{lem5}(i), $\tau_n\ge H$ with probability
tending to 1, hence by Lemma \ref{lem5}(ii), $\sup_{t\in[0,H]}|B_n(t)-t|\To0$.
Consequently, $\hat S_n\circ B_n\To X^{(2)}$, and so $X_n\To X$
in the uniform topology on $[0,H)$. Moreover, for $t\in[H,\tau_n)$, (on the event
$\tau_n>H$), $B_n(t)=B_n(H)+(t-H)$ by the non-idling property. It follows that
$(\hat S_n\circ B_n)(\cdot\w\tau_n)\To X^{(2)}(\cdot\w\bar\tau)$.
We now appeal to Skorohod's representation theorem,
by which we may assume without loss of generality that a.s.,
\begin{equation}\label{84}
\tau_n\to\bar\tau,\quad \text{ and }
(\hat\Xi_n,\hat\Xi_n\circ\bar E_n, (\hat S_n\circ B_n)(\cdot\w\tau_n))\to
(X^{(1)}, X^{(1)}(\tilde\mu\cdot), X^{(2)}(\cdot\w\bar\tau)),
\end{equation}
uniformly on compacts. If we let $X=X^{(1)}(\tilde\mu\cdot)-X^{(2)}$
then we also have $X_n(\cdot\w\tau_n)\to X(\cdot\w\bar\tau)$ a.s.

\begin{lemma}
  \label{lem6}
The estimate \eqref{82} holds with a sequence of RVs $\eps_n$ for which
$\E[\eps_n]\to0$ as $n\to\iy$.
\end{lemma}

\noi{\bf Proof of Theorem \ref{prop5}.}
The estimate \eqref{82} holds by Lemma \ref{lem6}. Hence \eqref{83} is valid.
Taking the limit inferior and using the convergence $\E[\eps_n]\to0$,
stated in Lemma \ref{lem6}, establishes the result.
\qed

\skp

We turn to the proofs of the lemmas.
We shall use uniform second moment bounds as follows. For every $t$,
\begin{equation}\label{85}
\sup_n\E[\|\hat\Xi_n\|_t^2]<\iy,
\qquad
\sup_n\E[\|\hat S_n\|_t^2]<\iy,
\end{equation}
where the first assertion follows by
Doob's $L^2$ maximum inequality \cite[\S 4.4]{Du10},
and the second is shown in \cite[Theorem 4]{KrTa1992}.

Recall that $c$ denotes a generic positive constant
(non-random, independent of $n$), whose value may change from
line to line. Moreover, $\{\Theta_n\}$ will denote a generic sequence of non-negative
RVs that are uniformly square integrable, i.e., $\sup_n\E[\Theta_n^2]<\iy$.
The value of the sequence $\{\Theta_n\}$ may also change from line to line.

Note that in view of \eqref{85}, for every fixed $t$, one has
$\|\hat S_n\|_t\le\Theta_n$ and $\|\hat\Xi_n\|_t\le\Theta_n$.
Moreover, $\|\hat S_n\circ B_n\|_t\le\Theta_n$
due to the fact $0\le B_n(s)\le s$ for all $s$.

\skp

\noi{\bf Proof of Lemma \ref{lem5}.}
We will argue that assertions (i)--(iii) hold for any sequence $\{T_i^n\}$
that achieves the limit inferior in \eqref{73}. Hence we fix such a sequence and denote it
by $\{T_i^n\}$ (thus, \eqref{80} is valid for this sequence).
On the other hand, the proof of part (iv) will require a certain construction;
it will be achieved by modifying this fixed sequence in a suitable way.

It follows from \eqref{79} and the fact that the sequence of costs is bounded
that
\begin{equation}
  \label{87}
  c_w\E[\hat W_n(1)]+c_w\E[\hat W_n(2)]+c_o\E[\hat\tau_n]\le c.
\end{equation}
We would like to deduce from \eqref{87}
that each of the terms on the LHS is bounded above by a constant.
Before we may do so we must provide a lower bound on each of these terms.
The first term is non-negative by its definition.

Next we show that $\E[(\hat\tau_n^-)^2]\le c$, equivalently $\hat\tau_n^-\le\Theta_n$.
By \eqref{61}, using the boundedness of $\bar E_n(H-)$
and $B_n(-)$, and the non-negativity of $I_n$, we have $\hat Q_n(H-)\ge-\Theta_n+p\hat E_n(H-)$.
By \eqref{60}, $\hat A^H_n\ge-\Theta_n-p\hat E_n(H-)$. Thus
\begin{equation}
  \label{88}
\hat Q_n(H)\ge-\Theta_n.
\end{equation}
Consider the event $\hat\tau_n\le 0$ and use \eqref{62}. On this event,
the expression $B_n(\tau_n)$ is bounded above by $\bar\tau$. Hence
the term $\hat S_n(B_n(\tau_n))$ is bounded in absolute value
by $\|\hat S_n\|_{\bar\tau}$.
Using this and the LB \eqref{88} in \eqref{62} gives
$\hat\tau_n1_{\{\hat\tau_n\le0\}}\ge-\Theta_n$. This gives
\begin{equation}
  \label{89}
  \E[(\hat\tau_n^-)^2]\le c, \qquad n\in\N.
\end{equation}
A lower bound on $\E[\hat W_n(2)]$ is achieved by considering the three expressions
$l_n=\E[\{\hat W_n(2)1_{\{\tau_n<H\}}]$, $l'_n=\E[\{\hat W_n(2)1_{\{\tau_n\in[H,\bar\tau]\}}]$,
$l''_n=\E[\{\hat W_n(2)1_{\{\tau_n>\bar\tau\}}]$.
For $l_n$, we use the lower bound $-cn^{1/2}$ on $\hat W_n(2)$, and \eqref{89}
by which $\PP(\tau_n<H)=\PP((\tau_n-\bar\tau)^->(\bar\tau-H))=\PP(\hat\tau_n^->n^{1/2}(\bar\tau-H))
\le cn^{-1}$. This shows $l_n\ge-cn^{-1/2}$.

For $l'_n$, on the event $\tau_n\in[H,\bar\tau]$, we use \eqref{64}.
By this equation, we have
\[
\hat W_n(2)1_{\{\tau_n\in[H,\bar\tau]\}}
\ge-\Theta_n-\frac{\mu}{2}n^{-1/2}(\hat\tau_n^-)^2.
\]
This shows $l'_n\ge-c$.

As for $l''_n$, it follows from \eqref{92} by a calculation similar to that leading
to \eqref{64} that, on $\{\tau_n>\bar\tau\}$,
\[
\hat W_n(2)
=\int_H^{\bar\tau}[q_n-\hat S_n(B_n(H)+t-H)+S_n(B_n(H))]dt
+n^{1/2}\int_{\bar\tau}^{\tau_n}Q_n(t)dt.
\]
The last term is non-negative on the indicated event, hence
\[
\hat W_n(2)1_{\{\tau_n>\bar\tau\}}\ge-\Theta_n.
\]
Thus $l''_n\ge-c$.
We conclude that $\E[\hat W_n(2)]\ge-c$.

In view of these lower bounds, \eqref{87} now implies
\begin{equation}\label{93}
(a)\ \E[\hat W_n(1)]\le c,
\qquad (b)\ \E[\hat W_n(2)]\le c,
\qquad (c)\ \E[\hat\tau_n^+]\le c.
\end{equation}
The bounds \eqref{89} and \eqref{93}(c) prove part (i) of the lemma.

Next, the tightness of $\hat\tau_n$, used in \eqref{62} implies the tightness
of $q_n$. By \eqref{94}, this gives the tightness of $n^{1/2}I_n(H-)$.
Since $I_n(H-)=I_n(H)=H-B_n(H)$, we obtain $B_n(H)\To H$. by the property
$|B_n(t)-B_n(s)|\le|t-s|$, this implies $\sup_{t\in[0,H]}|B_n(t)-t|\To0$,
and since for $t\in[H,\tau_n]$ we have $B_n(t)-B_n(H)=t-H$, the result stated
in part (ii) of the lemma follows.

The bound \eqref{93}(a) clearly implies the tightness of $\hat W_n(1)$.
By the expression \eqref{66} for $\hat W_n(1)$ and the convergence
$\tau_n\To\bar\tau$, this gives the tightness of $\int_0^H\hat Q_n(t)dt$.
Using \eqref{61+} and the Lipschitz property of $\Gamma_1$ (with constant 2),
\[
\int_0^H\Gamma_1[p\hat E_n](t)dt\le\int_0^H[p\hat E_n+X_n](t)dt+2H\|X_n\|_H
=\int_0^H\hat Q_n(t)dt+2H\|X_n\|_H,
\]
hence the tightness of the RHS implies that of the LHS. However,
the expression on the LHS is deterministic, thus it is, simply, bounded.
This gives the first assertion in part (iii).

The aforementioned tightness of the RVs $n^{1/2}I_n(H-)$,
along with the equality between the two expressions \eqref{94} and \eqref{95},
implies the tightness of the RVs $\Gamma_1[p\hat E_n(H-)]-p\hat E_n(H-)$
(in view of the tightness of the terms involving $\hat\Xi_n$ and $\hat S_n$
in these expressions). Arguing as above by the Lipschitz continuity
of $\Gamma_1$ establishes the second assertion in part (iii).

Finally we prove part (iv).
Given $\eps>0$, we first show that $\sup_{t\in[0,H-\eps]}(\bar E_n(t)-p^{-1}\mu t)\le\eps$
provided $n$ is sufficiently large.
Recall that $\bar E_n(t)$ is non-decreasing. Thus
if the inequality $\sup_{t\in[0,H-\eps]}(\bar E_n(t)-p^{-1}\mu t)>\eps$
is valid for some $n$ then there exists $t=t_n\in[0,H-\eps]$
such that for all $s\in[t,t+\eps_0\w\eps]$,
where $\eps_0=\eps p\mu^{-1}/2$,
$\bar E_n(s)-p^{-1}\mu s>\eps-p^{-1}\mu(s-t)>\frac{\eps}{2}$.
Hence by the definition of $\hat E_n$,
$\hat E^n(s)\ge \sqrt{n}\frac{\eps}{2}$ for the same set of times $s$.
By the definition of $\Gamma$, we have $\Gamma[p\hat E^n](s)\ge p\hat E_n(s)$,
for each $s$ in the interval alluded to above.
Hence $\int_0^H\Gamma(p\hat E^n)ds\ge c_\eps \sqrt{n}$, where $c_\eps>0$ depends on $\eps$
but not on $n$. By part (iii), this can occur for only finitely many $n$.
Hence the claim.

Next, given $\eps>0$, assume that the inequality
$\bar E_n(t)-p^{-1}\mu t>3\eps p^{-1}\mu$ is valid for some $n$ and $t=t_n\in[H-\eps,H)$.
Then $E_n(t)>p^{-1}\mu (t+3\eps)n\ge p^{-1}\mu (H+2\eps)n$.
Construct from $E_n$ another schedule
$E_n^{(1)}(s)=E_n(s)\w p^{-1}\mu (H+2\eps)n$, for $s\in[0,H)$.
Then the schedules agree on $s\in[0,t)$.
Moreover, the new schedule satisfies the constraint $\bar E_n(s)-p^{-1}\mu t<c\eps$ for all times
$t$.
The resulting queue length can only be decreased by this modification.
Moreover, the effect of the modification on the overage time is
negligible at the scaling limit, as follows by the following argument
which shows that under both schedules no idle time is accumulated during $[t,H)$,
with high probability. This owes to the fact that for $s\in[t,H)$,
\[
\bar E_n^{(1)}(s)=p^{-1}\mu(H+2\eps)
\]
hence according to \eqref{eq:queue-flln}, with the error term $\eps_n$
as in that equation,
\begin{align*}
\bar Q_n(s)&=p\bar E_n(s)-\mu s+\eps_n(s)+\sup_{u\leq s} (- p\bar E_n(u)+\mu u-\eps_n(u))^+\\
&= \mu(H+2\eps)-\mu s+\eps_n(s)+\sup_{u\leq s} (- p\bar E_n(u)+\mu u-\eps_n(u))^+\\
&\ge\mu(H+2\eps)-\mu s+\eps_n(s)\\
&\ge 2\mu\eps-\|\eps_n\|_H.
\end{align*}
This shows that for any $\eps>0$ one can construct $E_n$ for which
$\sup_{t\in[0,H)}(\bar E_n(t)-p^{-1}\mu t)<\eps$, while \eqref{80}
and assertions (i)--(iii) of the lemma hold. A diagonal argument
may now be used to take $\eps=\eps(n)\downarrow0$.

Next, a similar use of the formula \eqref{eq:queue-flln} now for the term
$\mu I_n$ shows that if
\[
\inf_{t\in[0,H]}(\bar E_n(t)-p^{-1}\mu t)\le -\eps
\]
then one has $I_n(H)\ge c\eps$ with high probability. In this case,
the cost associated with the scaled overage time
$\hat\tau_n$ grows without bound as $n\to\iy$.

Finally, the second assertion of part (iv)
follows from the first one by using the constraint $\bar E_n(H)=N_n$.
\qed

\skp

\noi{\bf Proof of Lemma \ref{lem6}.}
We first estimate $\hat W_n(1)$ from below.
By \eqref{61}, for $t\in[0,H)$, $\hat Q_n(t)$
is given by $p\hat E_n(t)+X_n(t)+n^{1/2}\mu I_n(t)$, where we recall that
\[
X_n(t)=\hat\Xi_n(\bar E_n(t))-\hat S_n(B_n(t)).
\]
We use \eqref{61+} and the fact that $\Gamma_1$ is Lipschitz with constant $2$
in the supremum norm to write
\begin{align*}
\hat W_n(1)&=\int_0^{H\w\tau_n}\hat Q_n(t)dt\\
&=\int_0^{H\w\tau_n}\{\Gamma_1[p\hat E_n+X_n](t)-\Gamma_1[p\hat E_n+X](t)\}dt
+\int_0^{H\w\tau_n}\Gamma_1[p\hat E_n+X](t)dt\\
&\ge
\int_0^H\Gamma_1[p\hat E_n+X](t)dt-\eps_n^0,
\end{align*}
where
\[
\eps_n^0=2(H\w\tau_n)\|X_n-X\|_{H\w\tau_n}
+\int_{H\w\tau_n}^H\Gamma_1[p\hat E_n+X](t)dt.
\]

The term $q_n$ (see \eqref{94}) appears both in the expression for $\hat W_n(2)$
and $\hat\tau_n$. We have
\begin{align}\label{86}
  q_n&=\Gamma_1[p\hat E_n+X_n](H-)+\hat\Xi_n(\al)
  -\hat\Xi_n(\bar E_n(H-))-p\hat E_n(H-)\\
  \notag
  &\ge\La_n - \eps_n^1,
\end{align}
where
\[
\La_n=\Gamma_1[p\hat E_n+X](H-)-p\hat E_n(H-)+X^{(1)}(\al)-X^{(1)}(\tilde\mu H),
\]
\[
\eps_n^1 = 2\|X_n-X\|_H+|\hat\Xi_n(\al)-X^{(1)}(\al)|
+|\hat\Xi_n(\bar E_n(H-))-X^{(1)}(\tilde\mu H)|.
\]
Recall that on the event $\tau_n<H$, $\hat W_n(2)=-n^{1/2}\bar W$.
On the event $\tau_n \ge H$, we have the expression \eqref{64}. We obtain (in both cases),
\begin{align*}
  \hat W_n(2)&\ge(\tau_n-H)^+(\hat Q_n(H-)+\hat A^H_n)
  -\int_H^{H\vee\tau_n}(X^{(2)}(t\w\bar\tau)-X^{(2)}(H))dt-\eps_n^2-\eps_n^3-\eps_n^4,
\end{align*}
where
\[
\eps_n^2=2\|\hat S_n\circ B_n(\cdot\w\tau_n)-X^{(2)}(\cdot\w\bar\tau)\|_{\tau_n},
\]
\[
\eps_n^3=\frac{\mu}{2}n^{1/2}(\tau_n-\bar\tau)^2,
\]
\[
\eps_n^4=n^{1/2}\bar W1_{\{\tau_n<H\}}.
\]
Hence
\[
\hat W_n(2)\ge(\bar\tau-H)\La_n-\int_H^{\bar\tau}(X^{(2)}(t\w\bar\tau)-X^{(2)}(H))dt
-\eps_n^2-\eps_n^3-\eps_n^4-\eps_n^5,
\]
where
\[
\eps_n^5=|\tau_n-\bar\tau|(\hat Q_n(H-)+\hat A_n^H)+\bar\tau\eps_n^1
+2|\tau_n-\bar\tau|\,\|X^{(2)}\|_{\bar\tau}.
\]
As for $\hat\tau_n$, using \eqref{62},
\[
\hat\tau_n\ge\mu^{-1}\Big[\La_n-X^{(2)}(\bar\tau)+X^{(2)}(H)-\eps_n^1-\eps_n^2\Big].
\]

We now combine the lower bounds obtained above on $\hat W_n(1)$, $\hat W_n(2)$
and $\hat\tau_n$ and compare to the expression \eqref{68} with $p\hat E_n$
substituted for $U$. The term $L(H-)$ appearing in \eqref{68} is related to
$\La_n$ via
\[
\E[L(H-)]=\E[\Gamma_1[p\hat E_n+X](H-)-p\hat E_n(H-)]=\E[\La_n].
\]
We have thus shown that \eqref{82} holds with
\[
\eps_n=c_5(\eps_n^0+\cdots+\eps_n^5),
\]
where $c_5$ is a constant.

It now remains to show that each of the terms $\E[\eps_n^i]$, $i=0,1,\ldots,5$
converges to zero.

The first term in $\eps_n^0$ converges to zero a.s.\ by \eqref{84}.
It follows from \eqref{85} that this term is uniformly integrable.
Hence its expectation also converges to zero.
As for the second term in $\eps_n^0$, using the Lipschitz
property of $\Gamma_1$, this term is bounded by
\[
1_{\{\tau_n<H\}}\Big\{\int_0^H\Gamma_1[p\hat E_n](t)dt+2H\|X\|_H\Big\}
\le1_{\{\tau_n<H\}}\Big\{c+2H\|X\|_H\Big\},
\]
where we used Lemma \ref{lem5}(iii) and $c$ is a suitable constant.
The indicator function converges to zero a.s., by \eqref{84},
while the expectation of $\|X\|_H$ is finite. Hence by dominated convergence,
$E[\eps_n^0]\to0$.

The a.s.\ convergence of $\eps_n^1$ and $\eps_n^2$
to zero follows from \eqref{84}
whereas their uniform integrability follows from \eqref{85}.
Thus $\E[\eps_n^1]\to0$ and $\E[\eps_n^2]\to0$.

Next, by \eqref{78}, $n^{1/2}(\tau_n-\bar\tau)^2=n^{-1/2}\hat\tau_n^2$.
Thus to show that $\E[\eps_n^3]\to0$ it suffices to improve the estimate
from Lemma \ref{lem5}(i) to show that $\E[\hat\tau_n^2]$ is bounded.

To this end, recall the expression \eqref{62}.
Since we already established the boundedness of the second moment of $q_n$,
it suffices to show that also the second and third
terms in \eqref{62} have bounded second moments.
Since $B_n(H)\le H$, the second moment
of the last term in \eqref{62} is bounded by $\E[\|\hat S_n\|_H^2]$;
hence \eqref{85} gives a uniform bound. It remains to show
that
\begin{equation}\label{81}
\sup_n\E[\hat S_n(B_n(\tau_n))^2]<\iy.
\end{equation}
Here we shall use the $3+\eps$ moment assumption.
First, note that $B_n(\tau_n)$
is the total time the server works on jobs,
thus is equal to the total arriving work. A bound on this
is given by $\sum_{i=1}^{N_n}\nu_{i,n}=n^{-1}\sum_{i=1}^{N_n}\nu_i$.
Denoting $w_n=\sum_{i=1}^{N_n}\nu_i$, we have
\begin{align*}
    \E[\hat S_n(B_n(\tau_n))^2]
    &=\E[1_{\{n^{-1}w_n<1\}}\hat S_n(B_n(\tau_n))^2]
    +\sum_{k=0}^\iy \E[1_{\{n^{-1}w_n\in[2^k,2^{k+1})\}}\hat S_n(B_n(\tau_n))^2]
    \\&\le
    \E[1_{\{n^{-1}w_n<1\}}\|\hat S_n\|^2_{n^{-1}w_n}]
    +\sum_{k=0}^\iy \E[1_{\{n^{-1}w_n\in[2^k,2^{k+1})\}}\|\hat S_n\|_{n^{-1}w_n}^2]\\
    &\le\E[\|\hat S_n\|^2_1]
    +\sum_{k=0}^\iy \E[1_{\{n^{-1}w_n\in[2^k,2^{k+1})\}}\|\hat S_n\|_{2^{k+1}}^2]\\
    &\le c+\sum_{k=0}^\iy \PP(w_n\ge n2^k)^{1/p}\,\E[\|\hat S_n\|_{2^{k+1}}^{2q}]^{1/q},
\end{align*}
where \eqref{85} is used for the first term,
and for the sum, H\"older's inequality is used, where $p^{-1}+q^{-1}=1$.
Fix $\beta\in(3,3+\eps)$. Then by Minkowski's inequality we have
$\E[w_n^\beta]\le cn^\beta\E[\nu_1^\beta]$. Hence $\PP(w_n\ge n2^k)\le c2^{-\beta k}$,
where $c$ is finite and does not depend on $n$ or $k$.
Next we appeal again to \cite[Theorem 4]{KrTa1992}, which states
that $\E[\|\hat S_n\|_t^\beta]^{1/\beta}\le c(t^{1/2}+1)$
(under the hypothesis that the $\beta$-th moment of $\nu_1$ is finite and $\beta\ge2$)
where $c$ does not depend on $t$ or $n$.
We use this estimate with $t=2^{k+1}$, and $q=\beta/2$ (accordingly, $p$
is determined). This gives
\begin{align*}
    \E[\hat S_n(B_n(\tau_n))^2]
    &\le c+c\sum_{k=0}^\iy 2^{-\beta k/p}(2^{k+1}+1).
\end{align*}
Now, $p=(1-2/\beta)^{-1}$, and since $\beta>3$, we have $p<3$.
In particular, $\beta>p$. Therefore the above sum is finite.
This proves \eqref{81}, hence follows the estimate on $\E[\eps_n^3]$.

To show that $\E[\eps_n^4]\to0$ amounts to showing that $n^{1/2}\PP(\tau_n<H)\to0$.
Now, $\tau_n\to\bar\tau>H$, and we have just shown that $\sup_n\E[\hat\tau_n^2]<\iy$.
Hence, for $c=\bar\tau-H$,
\[
n^{1/2}\PP(\tau_n<H)\le n^{1/2}\PP(n^{-1/2}|\hat\tau_n|>c)
\le \tilde cn^{-1/2},
\]
for some constant $\tilde c$. This shows $\E[\eps_n^4]\to0$.

Finally, for a bound on $\eps_n^5$ we use Cauchy-Schwartz to write
\[
\E[\eps_n^5]\le\E[(\tau_n-\bar\tau)^2]^{1/2}\E[q_n^2]^{1/2}
+\bar\tau\E[\eps_n^1]+c\E[(\tau_n-\bar\tau)^2]^{1/2}.
\]
We have already shown that $q_n$ are uniformly square integrable.
Hence the convergence of the first and last terms above
to zero follows from the boundedness of $\E[\hat\tau_n^2]$.
The second term has already been argued to converge to zero.

This concludes the proof that $\E[\eps_n^i]\to0$ for $i=0,1,\ldots,5$
and completes the proof of the lemma.
\qed

\subsection{Upper bound on the diffusion scale cost}\label{sec44}
In this section we propose a sequence of schedules (indexed by $n$)
whose diffusion scaled cost converges to the cost we obtained as a
solution to the BOP. This establishes that the BOP cost is an
asymptotic upper bound for diffusion scaled cost which, together with
the lower bound, establishes the asymptotic optimality of our proposed schedule.

Let $\beta^*$ be the optimal drift associated with the problem defined in Proposition \ref{prop4}, namely, $\beta^*$ is such that the control $U(t)=\beta^* t$ is optimal for the problem $\lim_{H\to\iy}\inf_{U\in\calU_H^{\rm lin}}\hat J_H(U)$. Note that $\beta^* < 0$. Let $n > (\beta^*/\mu)^2$, and for $t \in[0,H]$ define
\begin{equation} \label{eq:diff-opt-sched}
E^d_n(t) =
\left\{
\begin{array}{ll}
\left\lfloor\frac{nt}{p}\left(\mu+\frac{\beta^*}{\sqrt{n}}\right)\right\rfloor & t<H, \\
[0.2cm]
N_n & t=H.
\end{array}
\right.
\end{equation}
To see that $E^d_n(t)$ is an admissible schedule we need to verify that
$E_n^d(t)\geq 0$ and that it is non-decreasing in $t$. Both of these
requirements follow from the condition that $n>(\beta^*/\mu)^2$ and from
our assumption that $p^{-1}\mu H < \alpha$ in~\eqref{eq:overload}. The latter is used to
verify that the jump of $E_n^d$ at $t=H$ is non-negative.

Next, consider the diffusion-scaled schedule $\hat E^d_n := \sqrt n
\left( n^{-1} E^d_n - \l^{*}\right)$; recall that $\l^*$ is the fluid
optimal schedule, defined in~\eqref{eq:fop}. A straightforward computation shows
that, as $n \to \infty$
\begin{equation}
  \label{eq:diffusion-sched}
  \hat E^d_n(t) \to \hat u(t) :=
  \begin{cases}
    \beta^* t & t \in [0,H)\\
    0 & t \geq H.
  \end{cases}
\end{equation}
Furthermore, it can also be easily shown that the convergence is
uniformly on compact sets of the time index. It follows that as $H \to
\infty$, $\hat u$ converges to the drift of the stationary optimal RBM
associated with the problem $\lim_{H\to\iy}\inf_{U\in\calU_H^{\rm
    lin}}\hat J_H(U)$. The main result of this section proves that
this sequence of schedules asymptotically achieves the large time
horizon value of the BOP, determined in Proposition~\ref{prop4}.

\begin{theorem} \label{prop:diff-opt-sched}
  The sequence of schedules $\{T_{i,n}^d,~n\geq 1\}$ corresponding to
   the scheduling functions $\{E^d_n,~n\geq 1\}$ of \eqref{eq:diff-opt-sched}, satisfies
\[
\lim_{H \to \infty} \limsup_{n \to \infty} H^{-1} \hat
J_{n,H}(\{T_{i,n}^d\}) = \lim_{H\to \infty} \hat V_{H}.
\]
\end{theorem}

\noi{\bf Proof.}
Recall the diffusion-scaled cost in~\eqref{79}:
\[
\hat J_{n,H}(\{T_{i,n}^d\}) = \bbE[c_w \hat W_n(1) + c_w \hat W_n(2) + c_o
\hat \tau_n],
\]
where $\hat W_n(1), ~\hat W_n(2)$ and $\hat \tau_n$ are defined in
~\eqref{66},~\eqref{92} and~\eqref{68} (respectively). Applying
Skorokhod's representation theorem as in~\eqref{84} we have that
a.s.,
\[
\tau_n\to\bar\tau,\quad \text{ and }
(\hat\Xi_n,\hat\Xi_n\circ\bar E^d_n, (\hat S_n\circ B_n)(\cdot\w\tau_n))\to
(X^{(1)}, X^{(1)}(\tilde\mu\cdot), X^{(2)}(\cdot\w\bar\tau))~\text{as
} n \to\infty.
\]
As before, let $X := X^{(1)}(\tilde \mu \cdot) - X^{(2)}$, so that
$X_n(\cdot \wedge \tau_n) \to X(\cdot \wedge \bar \tau)$ a.s.

Consider $\hat W_n(1)$
first. From~\eqref{eq:overload},~\eqref{84}
and~\eqref{eq:diffusion-sched},
and the fact that $B_n(t) \in o(n^{-1/2})$ a.s. for all $t \in [0,H)$, it follows that $\hat W_n(1) = \int_0^{H\wedge \tau_n} \hat
Q_n(t) dt - n^{-1/2} B_n(\tau_n)$ satisfies:
\begin{equation}
  \label{eq:w1-limit}
  \hat W_n(1) \to \int_0^H \left( p \hat u + X \right) (t) dt
  ~\text{a.s. as } n \to\infty.
\end{equation}

Next, recall from~\eqref{64} that on the event $\{\tau_n  \geq H\}$
\[
\hat W_n(2) = \int_H^{\tau_n} \left[ q_n - \hat S_n(B_n(t)) + \hat
  S_n(B_n(H)) \right] dt - \frac{\m}{2} \sqrt n (\tau_n - \bar \tau)^2,
\]
where $q_n = \hat Q_n(H-) + \hat A_n^H$ and $\hat A_n^H = \hat
\Xi_n(\a) - \hat \Xi_n(\bar E^d_n(H-)) - p \hat E^d_n(H-)$. Again,
using~\eqref{84} and~\eqref{eq:diffusion-sched} it can be easily seen
that
\begin{equation}
  \label{eq:w2-limit-1}
  q_n \to \Gamma_1(p \hat u + X)(H-) + X^{(1)}(\a) - X^{(1)}(\tilde \mu
  H-) - p \hat u(H-) = \hat q~\text{a.s. as } n \to \infty.
\end{equation}
It follows from~\eqref{62},~\eqref{84} and~\eqref{eq:w2-limit-1} that
\begin{equation}
  \label{eq:w2-limit-2}
  \hat \tau_n \to \frac{1}{\m} \left( \hat q - X^{(2)}(\bar \tau) +
    X^{(2)}(H) \right) = \hat \tau ~\text{a.s. as } n \to \infty.
\end{equation}
Now, from the proof of Lemma~\ref{lem:upper-bound-convergence}(ii),
when traffic is scheduled per $E^d_n$ we have $\tau_n \to \bar \tau$
a.s. as $n \to \infty$. Furthermore, since $\bar \tau > H$, it follows
from~\eqref{64},~\eqref{eq:w2-limit-1} and~\eqref{eq:w2-limit-2} that
\begin{equation}
  \label{eq:w2-limit}
  \hat W_n(2) \to \int_H^{\bar \tau} \left[\hat q - X^{(2)}(t) +
    X^{(2)}(H) \right] dt~\text{a.s. as } n \to \infty,
\end{equation}
where we have used the fact that
\(
\sqrt n (\tau_n - \bar \tau)^2 = (\tau_n - \bar \tau) \hat \tau_n \to
0
\)
~a.s. as $n \to \infty$.

At this point, we have shown that, as $n \to \infty$, a.s.
\[
c_w \hat W_n(1) + c_w \hat W_n(2) + c_o \hat \tau_n \to c_w \int_0^H
\Gamma_1(p \hat u + X)(t) dt + c_w \int_H^{\bar \tau} \left[ \hat q -
  X^{(2)} (t) + X^{(2)}(H)  \right] dt + c_o \hat \tau.
\]
In order to prove convergence in $L^1$, following the analysis
in \S \ref{sec43} we will prove that the
second moments of $\hat W_n(1),~\hat W_n(2)$ and $\hat \tau_n$ are
bounded. This, however, follows directly from~\eqref{85} and the
analysis in the proof of Lemma~\ref{lem6}, and will not be repeated
here. Therefore, we have that
\[
\hat J_{n,H}(\{T_{i,n}^d\}) = c_w \hat W_n(1) + c_w \hat W_n(2)
+c_o \hat \tau_n
\]
is uniformly integrable, implying that $\lim_{n\to \infty} \hat J_{n,H}(\{T_{i,n}^d\}) =
  J_H(\hat u)$, where
\begin{align}
  \label{eq:L1-limit-cost}
  \notag
  J_H(\hat u ) &= \frac{c_w}{H} \bbE \left[ \int_0^H \Gamma_1\left( p \hat
    u + X\right)(t)  dt \right]\\
    &\quad + \frac{c_w}{H} \bbE \left[
  \int_H^{\bar \tau} \left(X^{(2)}(t) - X^{(2)} (H) + Q(H-) + \hat A^H
  \right) dt \right]
+ \frac{c_o}{H} \bbE\left[ \hat \tau \right].
\end{align}
Using the fact that $X^{(1)}$ and $X^{(2)}$ are Brownian motion
processes and the definition of $\hat \tau$ in~\eqref{eq:w2-limit-2}, $\hat J_H(\hat u)$ simplifies to
\[
\hat J_H (\hat u) = \frac{c_w}{H} \bbE \left[ \int_0^H \Gamma_1\left( p \hat
    u + X\right)(t)  dt \right] + \frac{\tilde c_o}{H} \bbE\left[ \Gamma_1(p
  \hat u + X)(H-) - p \hat u(H-) \right],
\]
where $\tilde c_o = c_o + c_w (\bar \tau - H)$.
Now, using the fact that
$\hat u$ converges to the drift function $\beta^* e$ as $H \to \infty$, it follows that $\lim_{H \to
  \infty} J_H(\hat u) = \lim_{H \to \infty} \hat V_H$, thus completing
the proof.
\qed

\subsection{The stochasticity gap at the diffusion scale}\label{sec45}
Continuing our investigation of quantifying the effect of the
stochasticity on the scheduling problem, we now
consider the asymptotic SG in the diffusion scale.
We study the limit of $\hat \gamma_{n,H} =
\sqrt n (V_{n,H} - V_{n,H}^{CI}) \geq 0$ as $n \to \infty$ and then $H
\to \infty$. We first show that the value of the CI problem in the
diffusion scale is asymptotically null.

\begin{lemma}~\label{lem:CI-diffusion-value}
  Fix $H > 0$. We have
\(
\hat \gamma_{n,H}^{CI} = \sqrt n (V_{n,H}^{CI} - \bar V_H) \to 0
\)
~as $n \to \infty$.
\end{lemma}
We delay the proof of the lemma to after the main result of this section.
\begin{proposition}~\label{thm:sg-diffusion}
Let Assumption \ref{assn1} hold. Then the
\textit{large horizon} SG is positive. More precisely,
  $\lim_{H\to \infty} \liminf_{n\to \infty} \hat \gamma_{n,H} =
  \lim_{H\to\infty} \hat V_H^* =: \hat
  V^*$.
\end{proposition}

\noi{\bf Proof.}
First note that $\hat \g_{n,H} = \sqrt n (V_{n,H} - \bar V_H) - \sqrt n (V_{n,H}^{CI} -
\bar V_H)$, where $\bar V_H$ is the FOP value from~\eqref{eq:FOP}. Also recall the definition $\hat V_{n,H} := \sqrt n (V_{n,H} - \bar
V_H)$. Then, Theorem~\ref{prop5} and Lemma~\ref{lem:CI-diffusion-value}
imply that
\begin{align}
  \label{eq:lbd}
  \lim_{H \to \infty} \liminf_{n \to\infty} \hat \g_{n,H} \geq \hat V^*.
\end{align}
On the other hand, Theorem~\ref{prop:diff-opt-sched} and
Lemma~\ref{lem:CI-diffusion-value} together imply that
\begin{align}
  \label{eq:ubd}
  \lim_{H \to \infty} \limsup_{n\to \infty} \hat\g_{n,H} \leq \hat V^*,
\end{align}
completing the proof.
\qed

\skp

\noi{\bf Proof of Lemma~\ref{lem:CI-diffusion-value}.}
Straightforward algebraic manipulation of~\eqref{57} shows that
\[
\begin{split}
\hat \g_{n,H}^{CI} = \bbE \bigg[ \frac{c_w}{H} \sqrt n \bigg(
    \frac{1}{n} \int_H^{H \vee \tau_n^*} \left( \Xi(N_n) - S_n(t) -1
    \right) dt &- \int_H^{\bar \tau} (\a p - \m t) dt \bigg) \\&  +
  \frac{c_o}{H} \sqrt n \left( (\tau_n^* - H)^+ - (\bar \tau - H)
  \right)\bigg],
\end{split}
\]
where we have used the fact that $\bar V_H = H^{-1}c_w \int_H^{\bar
  \tau} (\a p - \m t) dt+ H^{-1} c_o (\bar \tau - H)$.

In the proof of Lemma~\ref{lem:overage-CI} it was shown that $\tau_n^*
\to \bar \tau$ a.s. as $n \to \infty$. Let $x \in \bbR$ and consider the
event $\{\sqrt n (\tau_n^* - \bar \tau) > x\}$, or equivalently
$\{\tau_n^* > n^{-1/2} x + \bar \tau\}$. Recall that $\tau_n^* : =
\inf \{t > 0 : S_n(t) \geq \Xi(N_n) \}$, so that
\begin{eqnarray}
\nonumber
\{\tau_n^* > n^{-1/2} x + \bar \tau\} &=& \left \{ S_n\left( n^{-1/2} x
    + \bar \tau \right) < \Xi (N_n) \right\}\\
\label{eq:CI-event}
  &=& \{\hat \Xi_n(\a) - \hat S_n(n^{-1/2} x + \bar \tau) > x \m\},
\end{eqnarray}
where the last equality follows by simple algebraic manipulations and
recognizing that $\bar \tau = \m^{-1} \a p$.

Now, consider
\[
\hat \Xi_n(\a) - \hat S_n(n^{-1/2} x + \bar \tau) = \left( \hat \Xi_n(\a) -
\hat S_n(\bar \tau) \right) + \left(\hat S_n(\bar \tau) - \hat
S_n(n^{-1/2} x + \bar \tau)) \right.
\]
Recall that $(\hat \Xi_n, \hat
S_n) \Rightarrow (X^{(1)},X^{(2)})$ as $n \to \infty$. Therefore, it follows that $\hat \Xi_n(\a) -
\hat S_n(\bar \tau) \Rightarrow X^{(1)}(\a) - X^{(2)}(\bar \tau)$ as
$n \to \infty$. On
the other hand since $n^{-1/2} x + \bar \tau \to \bar \tau$ as $n \to
\infty$, we have $\hat S_n(\bar \tau) - \hat
S_n(n^{-1/2} x + \bar \tau) \To 0$ as $n \to \infty$. Therefore,
\begin{equation} \label{eq:CI-event-wc}
  \hat \Xi_n(\a) - \hat S_n(n^{-1/2} x + \bar \tau) \To X^{(1)}(\a) - X^{(2)}(\bar \tau).
\end{equation}
Displays~\eqref{eq:CI-event} and~\eqref{eq:CI-event-wc} together imply, as $n \to \infty$,
\begin{equation}
  \label{eq:CI-event-overage}
  \hat \tau_n^* \To \frac{1}{\m} \left( X^{(1)}(\a) - X^{(2)}(\bar \tau) \right).
\end{equation}

Now, consider $\hat{\underline{\tau}}_n^* = \sqrt n ((\tau_n^* - H)^+
- (\bar \tau - H))$. For each $n \geq 1$, we write
\(
\hat{\underline{\tau}}_n^* = \hat{\underline{\tau}}_n^*
1_{\{\tau_n^* > H\}} + \hat{\underline{\tau}}_n^* 1_{\{\tau_n^* \leq
  H\}}.
\)
The first term on the right hand side is simply $\hat \tau_n^*
1_{\{\tau_n^* > H\}}$. Since $\tau_n^* \To \bar \tau > H$ as $n \to
\infty$, it follows for large enough $n$ that $\tau_n^* > H$, implying
that $\hat{\underline{\tau}}_n^* \To \frac{1}{\m} \left( X^{(1)}(\a) -
  X^{(2)}(\bar \tau) \right)$. Following~\eqref{85}, it is
straightforward to deduce that $\hat{\underline{\tau}}_n^*$ is
uniformly integrable, implying that
\begin{align} \label{eq:5}
\bbE[\hat{\underline{\tau}}_n^*] \to \bbE\left[\m^{-1} \left( X^{(1)}(\a) -
  X^{(2)}(\bar \tau) \right)\right] = 0 ~\text{as}~n\to\infty.
\end{align}
Next, the first term in the definition of $\hat \g_{n,H}^{CI}$ can be
written as
\begin{flalign*}
 \begin{split}\sqrt n  \bigg(
    \frac{1}{n} \int_H^{H \vee \tau_n^*} & \bigg( \Xi(N_n) - S_n(t) -1
    \bigg) dt - \int_H^{\bar \tau} (\a p - \m t) dt \bigg) \\
 = \sqrt n &\bigg( \frac{1}{n} \int_H^{\tau_n^*} \left( \Xi(N_n)
      - S_n(t) \right) dt  - \frac{1}{n} (\tau_n^* - H) -
    \int_H^{\bar \tau} (\a p - \m t) dt \bigg) 1_{\{\tau_n^* > H\}}
  \\ & + \sqrt n \bigg( - \int_H^{\bar \tau} (\a p - \m t) dt \bigg) 1_{\{\tau_n^* \leq H\}}.
\end{split} &\\
\end{flalign*}
Consider the first term on the right hand side, under the event
$\{\tau_n^* > H\}$,
\begin{align*}
  \sqrt n \bigg( \frac{1}{n} \int_H^{\tau_n^*} \left( \Xi(N_n)
      - S_n(t) \right) dt  - \int_H^{\bar \tau} (\a p - \m t) dt
  \bigg) - \frac{1}{\sqrt n} (\tau_n^* - H),
\end{align*}
and focus on the term under the integral first:
\begin{align}\label{eq:CI-integral1}
  \sqrt n \left( \frac{1}{n} \int_H^{\bar \tau} \left( \Xi(N_n)
      - S_n(t) \right) dt  - \int_H^{\bar \tau} (\a p - \m t) dt
  \right) + \frac{1}{\sqrt n} \int_{\bar \tau}^{\tau_n^*} \left( \Xi(N_n)
      - S_n(t) \right) dt.
\end{align}
Note that the breaking up the integral is justified since $S_n(t)$ is
well-defined for all $t \geq 0$. Consider the latter integral:
\[
\frac{1}{\sqrt n} \int_{\bar \tau}^{\tau_n^*} \left( \Xi(N_n)
      - S_n(t) \right) dt = \int_{\bar \tau}^{\tau_n^*} (\hat
    \Xi_n(\a) - \hat S_n(t)) dt + \sqrt n \int_{\bar
      \tau}^{\tau_n^*}(\a p - \m t) dt.
\]
It follows that the right hand side equals
\[
\int_{\bar \tau}^{\tau_n^*} (\hat
    \Xi_n(\a) - \hat S_n(t)) dt + \a p \hat \tau_n^* - \frac{\m}{2}
    \hat \tau_n^* (\tau_n^* + \bar \tau),
\]
and as $n \to \infty$
\begin{align*}
  \frac{1}{\sqrt n} \int_{\bar \tau}^{\tau_n^*} \left( \Xi(N_n)
      - S_n(t) \right) dt \To \a p \hat \tau^* - \m \bar \tau \hat \tau^*.
\end{align*}
Now, using~\eqref{85} it can be shown that $\frac{1}{\sqrt n}
\int_{\bar \tau}^{\tau_n^*} \left( \Xi(N_n) - S_n(t) \right) dt$ is
uniformly integrable, so that
\begin{align}
  \label{eq:1}
  \bbE[\frac{1}{\sqrt n} \int_{\bar \tau}^{\tau_n^*} \left( \Xi(N_n)
      - S_n(t) \right) dt] \to \bbE[\a p \hat \tau^* - \m \bar \tau
  \hat \tau^*] = 0 ~\text{as}~n\to \infty.
\end{align}

Returning to~\eqref{eq:CI-integral1}, consider the term
\[
\sqrt n \left( \frac{1}{n} \int_H^{\bar \tau} \left( \Xi(N_n)
      - S_n(t) \right) dt  - \int_H^{\bar \tau} (\a p - \m t) dt
  \right)  = \sqrt n \left( \int_H^{\bar \tau} (\hat \Xi_n(\a) - \hat
    S_n(t)) dt \right).
\]
It is straightforward to see that this sequence of RVs converges
weakly to $\int_H^{\bar \tau} (X^{(1)}(\a) - X^{(2)}(t)) dt$. Furthermore, it is again true that the pre-limit integral is
uniformly integrable, implying that, as $n\to\infty$
\begin{align}
  \label{eq:4}
  \bbE \left[  \sqrt n \left( \int_H^{\bar \tau} (\hat \Xi_n(\a) - \hat
    S_n(t)) dt \right)\right] \to \bbE\left[ \int_H^{\bar \tau}
  (X^{(1)}(\a) - X^{(2)}(t)) dt \right] = 0.
\end{align}
From~\eqref{eq:5} and~\eqref{eq:4} it follows that $\hat \g_{n,H}^{CI}
\to 0$ as $n \to \infty$ for all $H > 0$.
\qed

\section{Conclusions}\label{sec5}
Exact solutions for the optimal scheduling problem studied in this
paper are intractable in general.
The analytical results provide the first rigorously
justified approximate solutions to this problem, in the large population
limit. We have taken the approach of first formulating
an optimization problem which is expected to govern the
asymptotics in the respective scales (based on formal limits),
then solving it, and finally proving that the value
of these limit problems indeed gives the limit of the values for the rescaled
scheduling problems. A byproduct of the last step, that is important in its own right,
is to derive asymptotically optimal schemes for the prelimit scheduling problems.

It is customary to distinguish a {\it control} problem,
where online information on the state of the system
is available to the decision maker, from an {\it optimization} problem,
where decisions are made at the initial time. A natural
decision/optimization model of the job scheduling problem is a
single-stage stochastic program without recourse \cite{ShDeRu2009}. The latter paradigm, however, assumes that
the distribution of the workload is fully specified. In general, it is much more natural to describe the arrival
and service processes and then infer the distribution of the workload process.
The three-step approach alluded to above is well established in work on
control in asymptotic regimes, specifically
in the heavy traffic literature, but is less studied in
optimization problems.
As far as the authors know,
an optimization problem involving diffusion, of the type
of the BOP we have formulated, has not been considered before in relation
to heavy traffic applications, or in the context of solving stochastic
programs without recourse. It seems that versions of this problem
might be relevant in applications far beyond the present
model. We also believe that our three-step
approach to solving the optimization problem may be used
to analyze and solve single-stage stochastic programs, where the
stochasticity is driven by functionals of multiple RVs.

While the FOP is easy to solve, we have not been able to find an explicit
solution to the BOP over a time interval of finite horizon.
As we have mentioned, the BOP is a convex optimization problem,
and thus it is plausible that one could treat it via numerical schemes;
this is left for future work.
However, one of our main findings is that, when set on an infinite time horizon,
this BOP is solvable explicitly. Its solution, in the form of a reflected
BM with constant drift, is rather simple. The form of the optimal
drift captures the tradeoff between the two parts of the cost.

Another main ingredient of our work
is the notion of a SG, that we have introduced as a means of
quantifying the performance loss due to the inherent stochasticity
in the model, as compared to the complete information problem.
As one may expect, we have shown that the gap converges to zero
in the fluid limit but remains positive in the diffusion limit.
It is natural to associate this to (but it certainly does not automatically follow from)
the fact that the FOP is a deterministic problem whereas the BOP
is stochastic. Moreover, our calculation of the gap in the diffusion limit
shows that it is proportional to the diffusion coefficient $\sigma$.
Thus the loss in performance due to stochasticity is proportional to
the standard deviation of the underlying noise.  It is also interesting to
note that the CI problem can be viewed as a single-stage stochastic
program \textit{with} recourse \cite{ShDeRu2009}, since the optimization is conducted
after the stochastic values are revealed to the decision-maker. The
SG thus provides a useful measure of the impact of
recourse on such problems.

A possible source of uncertainty not accounted for in this work,
but that is important in practice, is that of non-punctual arrivals.
This aspect may be addressed in future work.
Moreover, the analysis has focused exclusively on a
single server queue, and consequently the limit optimization problems
are one-dimensional.
Traffic scheduling in multi-queue networks are
natural to consider next, where the limit problems
are concerned with multidimensional diffusion processes
that are constrained to lie within a quadrant or, more generally, a cone.

\paragraph{Acknowledgment.}
The authors are grateful to Assaf Zeevi for the idea behind
the proof of Lemma \ref{lem4}.
The research of RA was supported in part by the Israel Science Foundation (grant 1184/16).
The research of HH was supported in part by the National Science Foundation (grant CMMI-1636069).

\footnotesize

\bibliographystyle{plain}
\bibliography{refs}

\begin{thebibliography}{10}

\bibitem{ArGl2012}
Victor~F. Araman and Peter~W. Glynn.
\newblock Fractional brownian motion with {$H < 1/2$} as a limit of scheduled
  traffic.
\newblock {\em J. Appl. Prob.}, 49(3):710--718, 2012.

\bibitem{BeJo2010}
Saif Benjaafar and Oualid Jouini.
\newblock Queueing systems with appointment-driven arrivals, non-punctual
  customers, and no-shows.
\newblock Technical report, University of Minnesota, 2011.

\bibitem{BG}
Dimitris Bertsimas and Vineet Goyal.
\newblock On the power of robust solutions in two-stage stochastic and adaptive
  optimization problems.
\newblock {\em Math. Oper. Res.}, 35(2):284--305, 2010.

\bibitem{Bi68}
Patrick Billingsley.
\newblock {\em Convergence of Probability Measures}.
\newblock Wiley, 1968.

\bibitem{Bo2007}
Vladimir~I. Bogachev.
\newblock {\em Measure Theory}, volume~1.
\newblock Springer, 2007.

\bibitem{CaVe2003}
Tugba Cayirli and Emre Veral.
\newblock Outpatient scheduling in health care: a review of literature.
\newblock {\em Prod. Oper. Mgmt.}, 12(4):519--549, 2003.

\bibitem{CaVeRo2006}
Tugba Cayirli, Emre Veral, and Harry Rosen.
\newblock Designing appointment scheduling systems for ambulatory care
  services.
\newblock {\em Healthcare Mgmt. Sci.}, 9(1):47--58, 2006.

\bibitem{ChYa01}
Hong Chen and David~D. Yao.
\newblock {\em Fundamentals of Queueing Networks: Performance, Asymptotics, and
  Optimization}.
\newblock Springer, 2001.

\bibitem{Du10}
Rick Durrett.
\newblock {\em Probability: Theory and Examples}.
\newblock Cambridge University Press, 2010.

\bibitem{EtKu2008}
Stewart~N. Ethier and Thomas~G. Kurtz.
\newblock {\em Markov Processes: Characterization and Convergence}, volume 282.
\newblock Wiley, 2009.

\bibitem{GuDe2008}
Diwakar Gupta and Brian Denton.
\newblock Appointment scheduling in health care: challenges and opportunities.
\newblock {\em IIE Tr.}, 40(9):800--819, 2008.

\bibitem{Ha2012}
Randolph~W. Hall.
\newblock {\em Handbook of Healthcare System Scheduling}.
\newblock Springer, 2012.

\bibitem{Ha1985}
J.~Michael Harrison.
\newblock {\em Brownian Motion and Stochastic Flow Systems}.
\newblock Wiley, 1985.

\bibitem{HaMe2008}
Refael Hassin and Sharon Mendel.
\newblock Scheduling arrivals to queues: A single-server model with no-shows.
\newblock {\em Mgmt. Sci.}, 54(3):565--572, 2008.

\bibitem{HoJaWa2015}
Harsha Honnappa, Rahul Jain, and Amy~R. Ward.
\newblock A queueing model with independent arrivals, and its fluid and
  diffusion limits.
\newblock {\em Queueing Syst.}, 80(1-2):71--103, 2015.

\bibitem{HoJaWa2016}
Harsha Honnappa, Rahul Jain, and Amy~R. Ward.
\newblock On transitory queueing.
\newblock arXiv:1412.2321, 2016.

\bibitem{KoFo1970}
Andrei~N. Kolmogorov and Sergei~V. Fomin.
\newblock {\em Introductory Real Analysis}.
\newblock Dover, 1970.

\bibitem{KrTa1992}
Elena~V. Krichagina and Michael~I. Taksar.
\newblock Diffusion approximation for {$GI/G/1$} controlled queues.
\newblock {\em Queueing Syst.}, 12(3):333--367, 1992.

\bibitem{LaLa2012}
Linda~R. LaGanga and Stephen~R. Lawrence.
\newblock Appointment overbooking in health care clinics to improve patient
  service and clinic performance.
\newblock {\em Prod. Oper. Mgmt.}, 21(5):874--888, 2012.

\bibitem{ShDeRu2009}
Alexander Shapiro, Darinka Dentcheva, and Andrzej Ruszczy{\'n}ski.
\newblock {\em Lectures on Stochastic Programming: Modeling and Theory}.
\newblock SIAM, 2009.

\bibitem{ZaPi2014}
Christos Zacharias and Michael Pinedo.
\newblock Appointment scheduling with no-shows and overbooking.
\newblock {\em Prod. Oper. Mgmt.}, 23(5):788--801, 2014.

\end{thebibliography}
\end{document}